\documentclass[letterpaper, paper, 10pt]{AAS}		

\usepackage{bm}
\usepackage{amsmath}
\usepackage{subfigure}
\usepackage{siunitx} 
\usepackage[colorlinks=true, pdfstartview=FitV, linkcolor=black, citecolor= black, urlcolor= black]{hyperref}
\usepackage{overcite}
\usepackage{footnpag}			      	
\usepackage{gensymb}
\DeclareMathAlphabet{\mathpzc}{OT1}{pzc}{m}{it}
\usepackage{mathtools}
\usepackage{adjustbox}
\usepackage{graphicx}
\usepackage{array}

\PaperNumber{19-424}

\begin{document}

\title{INTEGRATED OPTIMIZATION OF ASCENT TRAJECTORY AND SRM DESIGN OF MULTISTAGE LAUNCH VEHICLES}

\author{ \centering{Lorenzo Federici\thanks{PhD student, Department of Mechanical and Aerospace Engineering, ‘Sapienza’ University of Rome, Rome, Italy.},  
Alessandro Zavoli\thanks{Research Assistant, Department of Mechanical and Aerospace Engineering, ‘Sapienza’ University of Rome, Rome, Italy.},
Guido Colasurdo\thanks{Full Professor, Department of Mechanical and Aerospace Engineering, ‘Sapienza’ University of Rome, Rome, Italy.},}  \newline \centering
Lucandrea Mancini\thanks{Master student in Space Transportation Systems, ‘Sapienza’ University of Rome, Rome, Italy.}, \ and Agostino Neri\thanks{Head of VEGA-C Launcher System Integration and Launch System Synthesis Office, ESA ESRIN - VEGA Integrated Project Team, Frascati (Rome), Italy.}}

\maketitle{}

\begin{abstract}
This paper presents a methodology for the concurrent first-stage preliminary design and ascent trajectory optimization, with application to a \textit{Vega}-derived \textit{Light Launch Vehicle}. The reuse as first stage of an existing upper-stage (Zefiro 40) requires a propellant grain geometry redesign, in order to account for the mutated operating conditions. An optimization code based on the parallel running of several Differential Evolution algorithms is used to find the optimal internal pressure law during Z40 operation, together with the optimal thrust direction and other relevant flight parameters of the entire ascent trajectory. Payload injected into a target orbit is maximized, while respecting multiple design constraints, either involving the alone solid rocket motor or dependent on the actual flight trajectory. Numerical results for SSO injection are presented.
\end{abstract}

\section{Introduction}
The rapid evolution and innovation that the space industry has undergone in the last few years opened the possibility to reduce satellite dimensions and mass, especially for Low Earth Orbit (LEO) applications. 
As a main consequence of this new trend, several companies and space agencies have focused their attention on the development of a new family of launchers, here referred as \textit{Light Launch Vehicles} (LLVs), devoted to the orbital injection of light-weight payloads (\emph{i.e.}, with a total mass up to a few hundreds kilograms). The existence of these new simple and small launchers and the possibility of boarding multiple payloads within a single fairing, have permitted, and are going to permit even more in the future, to break down the cost of space launches.

The European launcher \textit{Vega} responds to the commercial market requirements for a new-generation light-weight launch vehicle capable of orbiting small to medium-size satellites.
The \textit{Vega} family is now going to be enlarged with the development of (at least) two new vehicles: \textit{Vega C} and \textit{Vega E}. The objective of such rockets is to increase the load capacity of \textit{Vega} and improve the overall performance, at more competitive costs.
At the same time, the aerospace community is thinking at ``scale versions'' of existing launch vehicles, devoted to smaller payload transportation into orbit. One notable example is the \textit{Vega C}-derived LLV here referred as \textit{Vega Light}, which  attempts at exploiting, with minimal changes, already existing and tested \textit{Vega} motors, in order to reduce development cost and increase reliability. 
\textit{Vega Light} is here considered to be a three solid stages launch vehicle composed of \textit{Zefiro} solid rocket motors (SRMs) of the class Z40, Z9 and ZX (X indicates a few tons of propellant), used as first, second and third stage, respectively. The attitude and orbital control system (AOCS) of the last stage has the duty of compensating SRM performance scattering and carrying out the final orbital injection. 
Z9 and ZX have been originally designed to operate as upper stages, thus they can be employed directly without any modification. Z40, instead, is now being developed as a second-stage motor for \textit{Vega C}: for this reason, a re-design of Z40 as a first stage (Z40FS) is mandatory to account for operation in mutated environmental conditions. This can be accomplished modifying as less as possible the present ground-tested SRM through an adjustment of the original nozzle and a modification of the propellant grain geometry. 

The initial geometrical configuration of the grain plays the main role in determining the performance of a SRM. Once the propellant type and nozzle geometry have been assigned, it is the only available means to achieve a suitable evolution of the burning surface, from which, in turn, the desired pressure and thrust history follow.
In the case of complex 3D propellant grains, such as finocyl grains, the design process involves (i) the mathematical modeling of the geometry and (ii) the following evaluation of various independent parameters that define the complex geometry. Changes in the value of each of these parameters bring significant effects on the thrust law, due to finocyl grain particular geometry.
In order to ensure that the best possible design (according to some reference performance index) against all achievable configurations is being acquired, the need arises for an intelligent parametric optimization that can control all the design variables.

The problem of optimizing the initial geometry of a three-dimensional finocyl grain has been faced by many investigations. Past works dealt with the problem by using a variety of optimization algorithms (sequential quadratic programming\cite{NisarConf2008}, simulated annealing\cite{Kamran2011}, hybrid optimization techniques\cite{Nisar2008, Adami2017} or hyper-heuristic approach\cite{Kamran2012}), but always assuming, as performance measure, a quantity referred to the SRM alone (such as average thrust\cite{NisarConf2008, Kamran2011}, specific impulse\cite{Adami2017}, propellant\cite{Nisar2008} or motor\cite{Kamran2012} gross mass) and considering constraints related only to SRM geometry or operation (burning time, fixed length or outer diameter of the grain, minimum thrust level, and so on). 
An attempt to the integrated design and optimization of a complete SRM-based aerospace system, such as an air-launched satellite vehicle\cite{Rafique2010}, a ground-launched interceptor\cite{Zeeshan2010} or a four-stage launcher\cite{Zafar2013}, which contemporarily considers solid propulsion, vehicle masses configuration and flight mechanics, was made, but by using simple models for both grain geometry and trajectory.
A concurrent motor robust design and launch vehicle trajectory optimization was carried out by using coupled direct/indirect\cite{Casalino2015} or meta-heuristic/indirect\cite{Casalino2014} optimization procedures, but, up to now, applied to a hybrid rocket engine upper stage, with 2D cylindrical geometries for the solid fuel.
The simultaneous optimization of the thrust profile and ascent trajectory was already applied to \textit{Vega} launcher\cite{Civek2017} and to a \textit{Vega}-derived LLV;\cite{Mancini2018} while the employed dynamical model and thrust law parameterization resemble what is presented here, substantial differences exist in the grain-design approach, flight strategy and optimization technique, because of the different aim of the present work. 

The main goal of this paper is, indeed, to present a new methodology for the concurrent first stage propellant grain preliminary design and complete ascent trajectory optimization of a multi-stage launch vehicle, which has the clear advantage to account for the close and mutual relationships that exist between trajectory and first-stage thrust law. 
In the present analysis, Z40 grain re-design will be accomplished by determining the optimal law for internal pressure during first stage operation, together with the optimal thrust direction along the whole trajectory. Relevant flight parameters (kick-off angle, coasting time-lengths, \emph{etc.}) are also optimized. The aim is to maximize the payload injected into a nominal target orbit, while respecting several path constraints, dependent on the first stage operation (\emph{e.g.}, maximum operating chamber pressure) and the actual trajectory (\emph{e.g.}, maximum dynamical-pressure), and terminal constraints (\emph{e.g.}, assigned orbital radius and inclination).
An optimization code based on the parallel running of several improved Differential Evolution algorithms is used to solve the constrained multidisciplinary optimization problem which arises from the formulation employed.
Numerical results are presented for a medium-altitude Sun-synchronous orbit mission.



\section{Trajectory Analysis}

\subsection{Dynamical Model}
A three-dimensional point-mass model (\textit{3-DoF} model) is considered for the launch vehicle. It has, indeed, the clear advantage to be time efficient with respect to a \textit{6-DoF} model, while allowing to capture the most prominent features of the mission. 
Under this assumption, at any time the rocket state is defined only by position $\bm{r}$, velocity $\bm{v}$ and mass $m$. 
The equations of motion, in an inertial reference frame, are:

\begin{align}
    \dot{\bm{r}} &= \bm{v} \\
    \dot{\bm{v}} &= \bm{g} + \frac{\bm{D} + \bm{T}}{m}\\
    \dot{m} &= - \dot{m}_{e}
\end{align}

A spherical model is assumed for the Earth; the gravitational acceleration vector is simply given by:
\begin{equation}
    \bm{g} = - \frac{\mu}{r^3} \bm{r}
    \label{eq:gravity}
\end{equation}
where $\mu$ is the Earth gravitational constant.

The expression for the drag force is, instead:
\begin{equation}
     \bm{D} = - \frac{1}{2}\, \rho(h)\, S_{ref}\, C_{D}(M) \, v_{rel}\, \bm{v}_{rel}
     \label{eq:Drag}
\end{equation}
where $S_{ref}$ denotes the vehicle cross-sectional area and $C_{D}(M)$ is the drag coefficient, considered here to be just a function of the Mach number $M$. 


A simplified atmospheric model is considered; air density $\rho$, pressure $p_a$ and temperature $T_a$ are evaluated as a function of the altitude $h$, according to U.S. Standard Atmosphere 1976 model.\footnote{\url{https://ntrs.nasa.gov/archive/nasa/casi.ntrs.nasa.gov/19770009539.pdf}}
Moreover, it is assumed that the atmosphere rigidly rotates with the Earth, with spin rate $\bm{\omega}_E$. Thus, the vehicle velocity relative to the atmosphere is calculated as: $\bm{v}_{rel} = \bm{v} - \bm{\omega}_{E} \times \bm{r}$.
The lift force is here neglected, as it is usually done during the preliminary \textit{3-DoF} trajectory analysis.

Concerning the thrust force generated by SRMs, it can be expressed as \cite{Wiesel1997}:
\begin{equation}
    \bm{T} = \left [ \dot{m}_{e} u_{e} + (p_{e} - p_{a}) A_{e} \right ] \bm{\hat{T}} = \left ({T}_{vac} - p_{a} A_{e} \right )  \bm{\hat{T}}
    \label{eq:T}
\end{equation}
where $u_e$ indicates the combustion products exhaust velocity, with a mass flow rate equal to $\dot{m}_{e}$, $A_e$ the nozzle exit area, and $p_e$ the ejection pressure. $T_{vac}$ denotes, instead, the thrust in vacuum.

The differential equations and all the involved constants and variables are made dimensionless; in this way, relevant quantities are in a small range around unity, reducing the errors related to the finite digit arithmetic of the computer.

\subsection{Path Constraints}
In order to ensure the structural integrity of both vehicle and payload during the entire ascent trajectory, a variety of constraints must be set during the simulation phase. Such constraints are, in general, non-linear functions of the values that the state variables assume, and their main effect is the reduction of the space of the viable trajectories.

In particular, the following path constraints must be respected along the atmospheric flight:
\begin{align}
q = \frac{1}{2} \rho v^2_{rel} &\leq q^{max} & \mbox{Dynamic pressure} \label{eq:Con1}\\
(q \ast \alpha) = q \alpha &\leq (q \ast \alpha)^{max} & \mbox{Bending load} \label{eq:Con2}\\
\Dot{Q}_W = \frac{1}{2} \rho v^3_{rel} &\leq \Dot{Q}_W^{max} & \mbox{Heat flux after fairing jettisoning} \label{eq:Con3}\\
a_x = \frac{\dot{\bm{v}} \cdot \bm{\hat{T}}}{ g_{0}} &\leq a^{max}_{x} & \mbox{Axial acceleration} \label{eq:Con4}  
\end{align}

By assuming, for simplicity, that the thrust is always aligned with rocket symmetry axis, the angle of incidence can be evaluated as:
\begin{equation}
        \alpha = \arccos{(\bm{\hat{T}} \cdot \bm{\hat{v}}_{rel})}
        \label{eq:incidence}
\end{equation}

Actually, the launcher and the payload are subjected, during the flight, to both static and dynamic loads.
For this reason, the so-called \textit{quasi-static load} (QSL), \emph{i.e.}, the most severe combination of dynamic and static axial accelerations encountered during the mission, might be considered as axial load.\cite{Vega2012}
The dynamic loads acting of the vehicle are typically higher during the flight of the first stage, because of the greater thrust level and the interaction with the atmosphere. For this reason, two different values $a^{max, 1}_{x}$ and $a^{max, 2}_{x}$ for the maximum reachable level of axial acceleration have been considered during the first-stage flight and the rest of the ascent, respectively.
The maximum admissible values assumed for the variables in Eqs. from \eqref{eq:Con1} to \eqref{eq:Con4} are reported in Table \ref{tab:Constraints}.

\begin{table}[h]
   \caption{Constraints assumed for \textit{Vega Light}}
    \label{tab:Constraints}
    \centering
    \begin{tabular}{c|c|c|c|c}
      \hline 
       $q^{max}$ & $(q \ast \alpha)^{max}$ & $a^{max, 1}_{x}$ & $a^{max, 2}_{x}$ & $\dot{Q}_W^{max}$\\
        $[Pa]$ & $[Pa \cdot \degree]$ & $[g]$ & $[g]$ & $[W/m^2]$\\
      \hline 
      54000 & 78000 & 4 & 5 & 900\\
      \hline 
    \end{tabular}
\end{table}


\subsection{Terminal Constraints} 
In order to ensure the payload injection into a circular orbit of given radius $\tilde{r}$ and inclination $\tilde{i}$, 
the following terminal constraints, acting on final periapsis radius $r_{p}$, apoapsis radius $r_a$ and inclination $i$, must be considered:
\begin{equation}
\begin{array}{lcl}
    r_{p} & \geq & \tilde{r} -  \delta  r\\
    r_{a} & \leq & \tilde{r} +  \delta  r\\
    |i - \tilde{i}| & \leq & \delta  i
    \end{array}
    \label{eq:Terminal}
\end{equation}
In the present analysis, an injection accuracy of $\delta r = 50\,m$ on radius, and $\delta  i = 0.01\,\degree$ on inclination are enforced; the right ascension of the ascending node, instead, can be easily adjusted to the desired value by properly choosing the launch time.

\subsection{Flight Strategy} 
\begin{figure} [h!]
    \centering
    \includegraphics[width=0.4\textwidth]{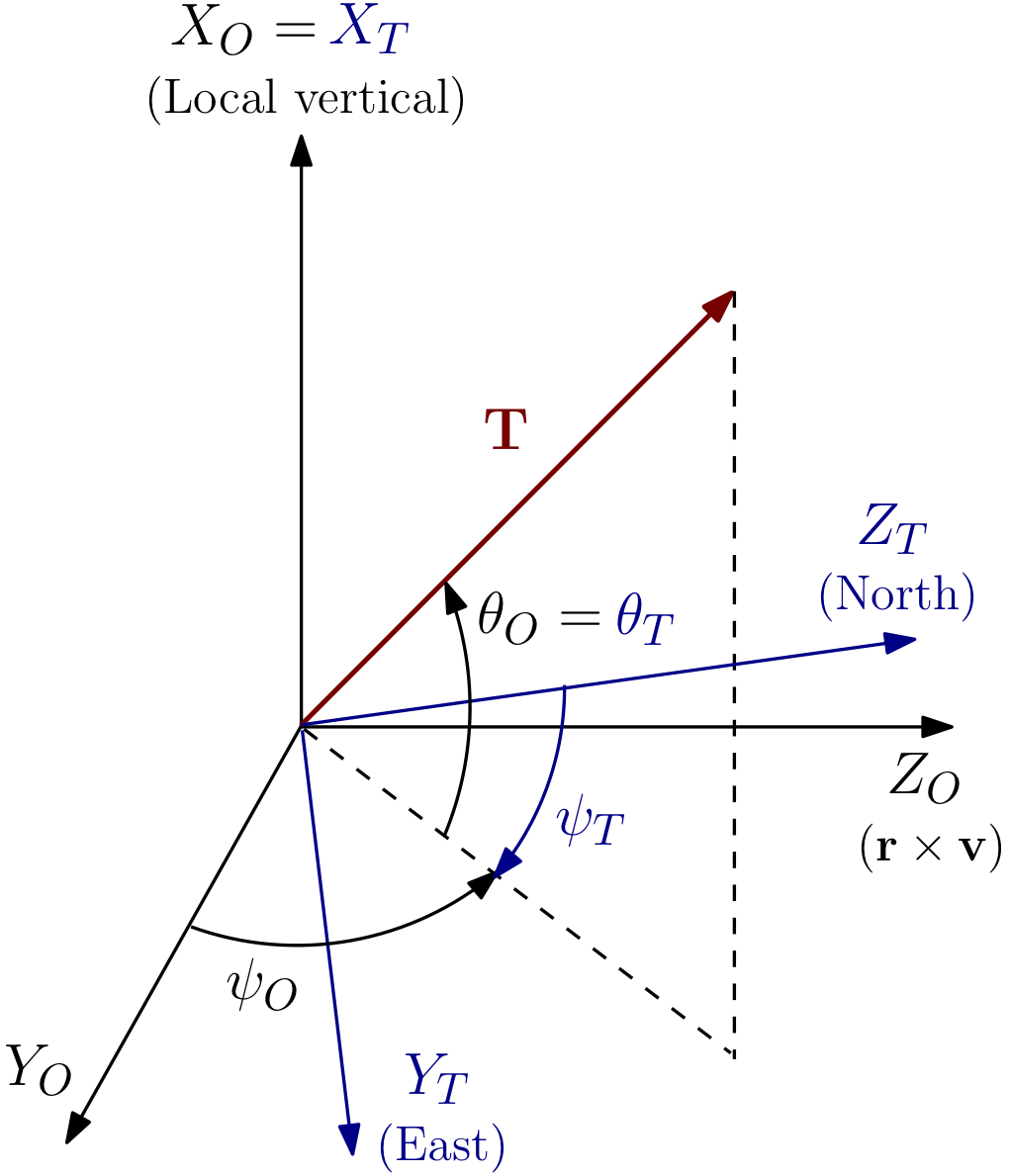}
    \caption{Thrust $\bm{T}$ in terms of flight path angle $\bm{\theta}$ and azimuth $\bm{\psi}$ in topocentric $\bm{\mathcal{F}_T}$ and RTN $\bm{\mathcal{F}_{O}}$ frames.}
    \label{fig:frames}
\end{figure}

In the case of an aerospace vehicle, devising a flight strategy or guidance program means to schedule the instantaneous direction of the generated thrust $\bm{\hat{T}}$ with respect to a reference frame along the entire trajectory.
As guidance variables defining $\bm{\hat{T}}$, the thrust flight path angle $\theta$ and azimuth $\psi$ have here been adopted. Such variables will be referred either to the topocentric frame $\mathcal{F}_T$ (topocentric flight path angle $\theta_T$ and azimuth $\psi_T$) or to the Radial-Transverse-Normal (RTN, sometimes defined ``orbital'') frame $\mathcal{F}_{O}$ (orbital flight path angle $\theta_{O}$ and azimuth $\psi_{O}$), according to convenience. Both reference frames are centered in the instantaneous position of the rocket; they are shown, together with the above-mentioned angles, in Figure~\ref{fig:frames}.

The adopted reference guidance program is roughly the same used by the majority of the launchers currently on the market, and, in particular, by launcher \textit{Vega}.\cite{Vega2012} It is used, as well, with light modifications, in many ascent trajectory optimization problems.\cite{Markl2001, Pagano2010, Brusch1976, Cremaschi2012}. 
This program consists of a series of specific guidance laws, each to be performed in a determined leg of the trajectory; so, it can also be seen as a way to split the ascent trajectory into phases or \textit{arcs}.
The following guidance laws make up the selected guidance scheme for \textit{Vega Light} LV.

\subsubsection{Lift-Off.} 
The natural way to \textit{lift-off} is, of course, to fly a vertical trajectory; the topocentric flight path angle is fixed to $90\degree$ while, because of the vertical direction of the thrust, the topocentric azimuth is not defined. A vertical \textit{lift-off} is mandatory since, at the main engine ignition, the generated thrust is the only force able to fight gravity and detach the vehicle from the ground. However, since the \textit{lift-off} phase has high gravity and misalignment losses, it is flown for the minimum time necessary to reach a safe height above the ramp. Then, the \textit{lift-off} time-length $t_{1,1} = t_1 - t_0$ is fixed to few seconds.

\subsubsection{Pitch-Over and Recovery.} 
Right after \textit{lift-off}, the vehicle starts to rotate by deflecting the engine nozzle. As a consequence, the pitch angle (equal to thrust flight path angle in a \textit{3-DoF} model) decreases (\textit{pitch-over} phase).
A simple linear \textit{pitch-over} law is used:
\begin{equation}
\begin{array}{lcl}
    \theta_{T} = 90\degree + \frac{t - t_1}{t_2 - t_1} \left(\theta_{T}(t_2) - 90\degree \right) & \mbox{for} & t_1 \leq t \leq t_2 
\end{array}
\label{eq:POlaw_FPA}
\end{equation}
The \textit{kick angle} $\theta_{T}(t_2)$ and the \textit{pitch-over} time-length $t_{1,2} = t_2 - t_1$ are optimized. The range in which the optimizer looks for them is reported in Table~\ref{tab:GuidLaws}. A constraint is set on the maximum admissible value for the \textit{pitch-over} rate, because of the limited maneuverability of the engine nozzle:
\begin{equation}
\dot{\theta}_{PO} = \frac{90\degree - \theta_T(t_2)}{t_{1,2}} \leq \dot{\theta}^{max}_{PO}
\label{eq:ConPO}
\end{equation}
with $\dot{\theta}^{max}_{PO} = 2\, \degree/s$. It will be considered as an additional path constraint.

Concerning the azimuth, a constant value is used in $\mathcal{F}_T$. The optimizer will look for it in a range centered in $\psi^{ref}_{T}$, \emph{i.e.}, the optimal topocentric azimuth in the case of non-rotating Earth:
\begin{equation}
    \psi^{ref}_{T} = \pi - \arcsin{\left(\frac{\cos{\tilde{i}}}{\cos{\delta_{LB}}}\right)}
    \label{eq:POlaw_Az}
\end{equation}
with $\delta_{LB} = 36.97\degree N$ the latitude of the launch base, supposed to be in the Azores islands.

A \textit{recovery} phase is here considered at the end of the \textit{pitch-over} maneuver, in order to progressively align thrust and relative velocity before the gravity turn starts. The thrust direction is fixed in $\mathcal{F}_T$, while the best \textit{recovery} time-length $t_{1,3} = t_3 - t_2$ is determined by the optimization process.

\subsubsection{Gravity Turn and Coasting.} 
The \textit{Zero Lift Gravity Turn} (ZLGT) or, simply, \textit{gravity turn}, is a clever maneuver, widely exploited by launchers, during which the thrust vector $\bm{T}$, the roll axis and the relative velocity $\bm{v}_{rel}$ are kept aligned to preclude the generation of any side force, and the yaw component of gravity is exploited to rotate slowly the velocity vector downward.\cite{Wiesel1997}
The ZLGT guidance law is:
\begin{equation}
    \bm{\hat{T}} = \bm{\hat{v}}_{rel} 
\label{eq:ZLGTlaw}
\end{equation}

The first stage ZLGT (\textit{ZLGT 1st stage}) is performed up to time $t_4 = t_0 + t_{B,1}$, where $t_{B,i}$ denotes the burning time of the \emph{i-th} stage. Because of the first stage grain redesign, the time $t_{B,1}$ is part of the optimization variables. 
At this point, the engine shut-down occurs, followed by the first stage separation from the rest of the rocket. 
A short \textit{coasting} phase (\textit{coasting 1-2}) with fixed duration $t_{C1} = t_5 - t_4$ is considered to allow the separation before the second stage ignition. Along this arc, the vehicle performs a ballistic flight, with null thrust.
A ZLGT trajectory is flown by the second stage too (\textit{ZLGT 2nd stage}), for a total time equal to its burning time $t_{B,2} = t_6 - t_5$ (now fixed).
When the second stage ends its propellant, it separates from the rest of the rocket. At this time, also the fairing jettisoning occurs. The following ballistic \textit{coasting} phase (\textit{coasting 2-3}) takes place close to the atmosphere limit ($h \approx 150\,km$) and it is characterized by a high relative velocity ($v_{rel} \approx 4\,km/s$). Therefore, it can be exploited to make the vehicle gain altitude, at the expense of a kinetic energy decrease.
The \textit{coasting 2-3} duration $t_{C2} = t_7 - t_6$ is determined by the optimization algorithm.

\subsubsection{Third Stage Flight.}
\begin{figure} [h!]
    \centering
    \includegraphics[width = 0.5\textwidth]{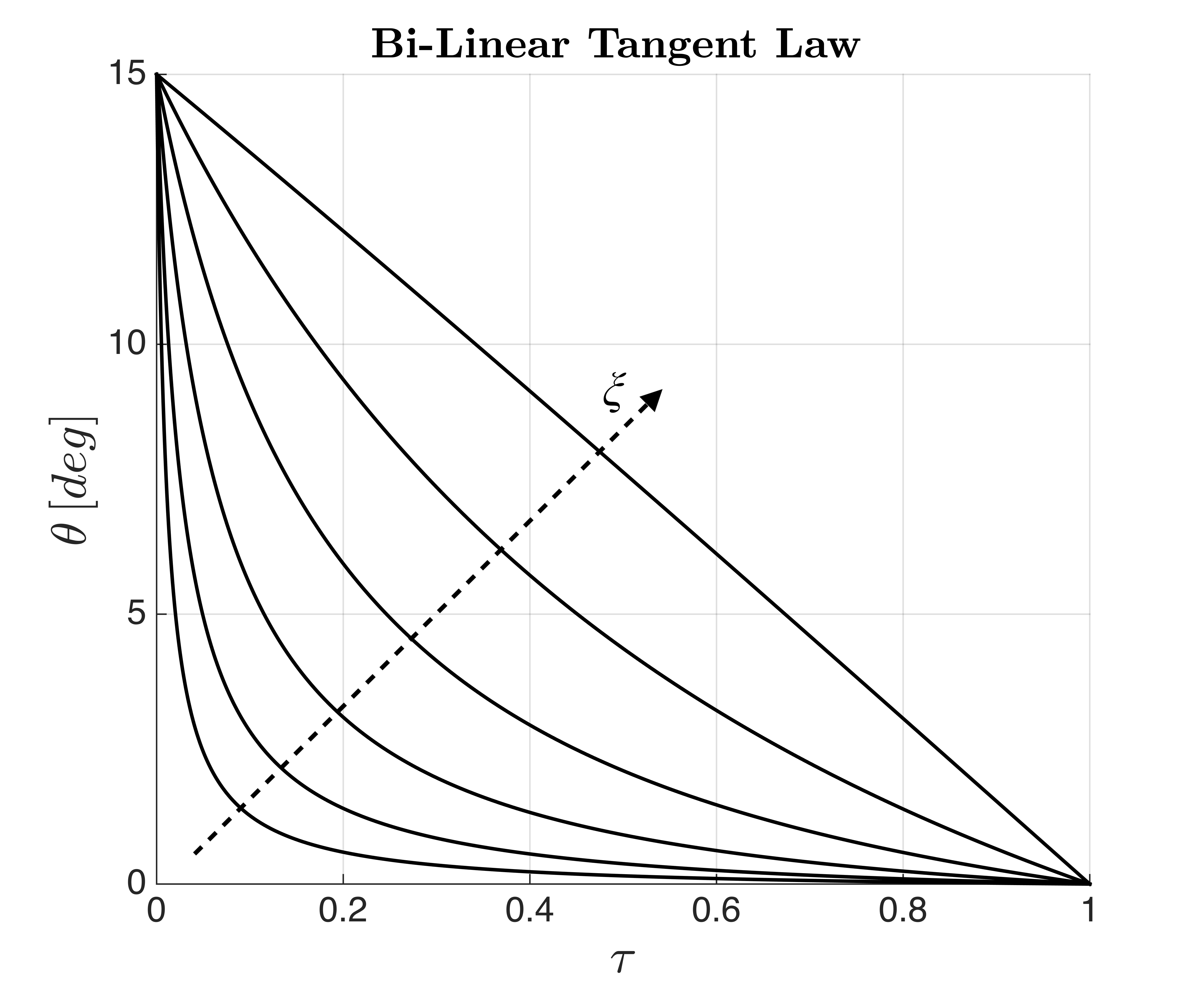}
    \caption{Effect of the curvature parameter $\bm{\xi}$ with base $\bm{a = 100}$ on the BLT. $\bm{\xi}$ values are sampled with a $\bm{0.2}$ interval.}
    \label{fig:BLT}
\end{figure}

The third stage ignition marks the end of the atmospheric flight of the vehicle.
A guidance law for the flight path angle which can be used, with good results, during the extra-atmospheric flight, is the so-called \textit{Bi-Linear Tangent law} (BLT).\cite{Markl2001,Brusch1979}
Strictly speaking, the BLT is the optimal guidance law under the simplifying hypothesis of flat Earth, plane motion and no aerodynamic forces: 
\begin{equation}
    \tan{(\theta)} = \frac{C_2 t + C_3}{C_1 t + 1}
\end{equation}
A good compromise is to express the value of the three constants $C_1,\,C_2,\,C_3$ in terms of the initial and final flight path angles and a measure of the curvature $c = a^\xi$, with $a = 100$ and $\xi \in [-1,\, 0]$.\cite{Markl2001}
Therefore, the expression of the guidance law for the flight path angle is, in RTN frame:
\begin{equation}
    \tan{ \left (\theta_{O} \right )}  =  \frac{a^\xi \tan{\left (\theta_{O}(t_7)\right )} + \left [\tan{\left (\theta_{O}(t_8)\right )} - a^\xi \tan{\left (\theta_{O}(t_7)\right )} \right] \tau}{a^\xi + \left (1 - a^\xi\right ) \tau}
\label{eq:BTLlaw_FPA}
\end{equation}
where: $\tau = \frac{t - t_7}{t_8 - t_7} \in [0,\, 1]$ is the non-dimensional time in the chosen interval.
Figure~\ref{fig:BLT} shows the trend of the BLT law for different values of $\xi$, by supposing $\theta_{O}(t_7) = 15\degree$ and $\theta_{O}(t_8) = 0\degree$. For $\xi = 0$, the \textit{linear tangent law} will result.
The optimized variables are reported in Table~\ref{tab:GuidLaws}.

Concerning the azimuth, a constant value is assumed in the RTN frame, and determined through the optimization process.
The time-of-flight of the third stage is imposed equal to its (fixed) burning time: $t_{B,3} = t_8 - t_7$.

\subsubsection{Orbital Injection.}
\label{sec:OrbInj}
The third stage of \textit{Vega Light} remains connected to the payload until its final \textit{orbital injection}.
A third ballistic \textit{coasting} (\textit{coasting 3-4}) lets the vehicle gain altitude before the final burn is performed. Its duration $t_{C4} = t_9 - t_8$ could vary from few seconds to several minutes, according to the altitude of the target orbit.
The final \textit{orbital injection} of the payload is carried out by means of a single shoot of the AOCS of ZX motor.
A simple constant thrust direction in RTN frame, determined by the optimization process, is selected for this final arc.
The burning time of the AOCS $t_{B,4} = t_{10} - t_9$ is evaluated so that its propellant is entirely consumed.
A summary of the described guidance laws can be seen in Table~\ref{tab:GuidLaws}.

\begin{table}[h!]
  \centering
   \caption{Guidance program for \textit{Vega Light} trajectory.} 
    \label{tab:GuidLaws}
    \begin{tabular}{c|c|c|c|c}
     \hline
       \textbf{\textit{Arc}} & \textbf{Phase} & \textbf{Time} & \textbf{Optimization} & \textbf{Variables}\\
        & & \textbf{interval} & \textbf{variables} & \textbf{boundaries}\\
      \hline
       $\bf 1$ & \textbf{Lift-off} & $[t_0,\, t_1]$  & none & none\\
       \hline
       $\bf 2$ &  \textbf{Pitch-over} & $[t_1,\, t_2]$ & $t_{1,2}$ & $[5\,s,\, 15\,s]$\\
       &  & & $\theta_T(t_2)$ & $[70\degree,\, 89.9\degree]$\\
       &  & & $\psi_T(t_2)$ & $[\psi^{ref}_T - 5\degree,\,\psi^{ref}_T + 5\degree] $\\
       \hline
       $\bf 3$ &  \textbf{Recovery} & $[t_2,\, t_3]$ & $t_{1,3}$ & $[1\,s,\, 8\,s]$\\
       \hline
       $\bf 4$ &  \textbf{ZLGT 1st stage} & $[t_3,\, t_4]$ & $t_{B,1}$ & $[85\,s,\, 120\,s]$\\
       \hline
       $\bf 5$ &  \textbf{Coasting 1-2} & $[t_4,\, t_5]$ & none & none\\
       \hline
       $\bf 6$ &  \textbf{ZLGT 2nd stage} & $[t_5,\, t_6]$ & none & none\\
       \hline
       $\bf 7$ &  \textbf{Coasting 2-3} & $[t_6,\, t_7]$ & $t_{C2}$ & $[1\,s,\, 100\,s]$\\
       \hline
       $\bf 8$ &  \textbf{3rd stage} & $[t_7,\, t_8]$  & $\theta_O(t_7)$ & $[-15\degree,\,15\degree]$\\
       & \textbf{flight} &  & $\theta_O(t_8)$  & $[-15\degree,\,15\degree]$\\
       &  &  & $\xi$  & $[-1,\,0]$\\
       &  & & $\psi_O(t_7)$  & $[-180\degree,\,180\degree]$\\
       \hline
       $\bf 9$ &  \textbf{Coasting 3-4} & $[t_8,\, t_9]$ & $t_{C4}$ & $[1\,s, 1.5\,h]$\\
       \hline
       $\bf 10$ &  \textbf{Orbital} & $[t_9,\, t_{10}]$ & $\theta_O(t_9)$ & $[-10\degree,\,10\degree]$\\
       & \textbf{injection} &  & $\psi_O(t_9)$  & $[-180\degree,\,180\degree]$\\
      \hline
    \end{tabular}
\end{table}

\section{First Stage SRM Grain Design} 

\subsection{Solid Propulsion Model} 
A simplified 1-D stationary isentropic model for the expansion process in the nozzle is adopted to evaluate Z40FS performance, \emph{i.e.}, mass flow rate $\dot{m}_e$ and vacuum thrust $T_{vac}$, as a function of equilibrium chamber pressure $p_c$ and temperature $T_c$, combustion products properties (molecular mass $\mathcal{M}$, specific heat ratio $\gamma$) and nozzle geometry (throat area $A_t$ and exit area $A_e$), as it is typically done during a preliminary phase\cite{Cornelisse1979}. 

Under such hypothesis, the expressions for vacuum thrust and mass flow rate are:
\begin{equation}
       T_{vac}  = \eta_{thr}\, \Gamma p_c A_t \sqrt{\left (\frac{2 \gamma}{\gamma - 1}\right)\, \left [1 - \left (\frac{p_e}{p_c}\right)^\frac{\gamma - 1}{\gamma}\right]} + p_e A_e = p_c A_t C_{F, vac}
\label{eq:THRvac}
\end{equation}
\begin{equation}
    \dot{m}_e = \frac{\Gamma}{\eta_{c^\ast}} \frac{p_c A_t}{\sqrt{\frac{\mathcal{R}}{M} T_c}} = \frac{p_c A_t}{c^\ast}
\label{eq:MFR}
\end{equation}
with $\mathcal{R}$ the universal gas constant and $\Gamma = \sqrt{\gamma} \left (\frac{2}{\gamma + 1}\right )^\frac{\gamma + 1}{2(\gamma - 1)}$. $C_{F, vac}$ and $c^\ast$ denote the vacuum thrust coefficient and the characteristic velocity, respectively.
The quality factors $\eta_{thr}$ and $\eta_{c^\ast}$ have been introduced to account for losses in performance of the motor due to its non-ideal behaviour.
Pressure ratio $p_e/p_c$ between the exit and entry section of the nozzle is a function of propellant type and nozzle expansion ratio $\epsilon = A_e/A_t$ through the relation:
\begin{equation}
    \frac{A_e}{A_t} = \frac{\Gamma}{\sqrt{\left (\frac{2 \gamma}{\gamma - 1}\right)\, \left ( \frac{p_e}{p_c}\right )^\frac{2}{\gamma} \left [1 - \left (\frac{p_e}{p_c}\right)^\frac{\gamma - 1}{\gamma}\right]}} 
\label{eq:expansion}
\end{equation}
The vacuum specific impulse is, instead, by definition: $I_{sp, vac} = \frac{T_{vac}}{\dot{m}_e g_0}$, with $g_0$ the reference sea-level gravitational acceleration.
A quasi-steady zero-dimensional internal ballistic model is here adopted for predicting the physical quantities behaviour within the combustion chamber: gas properties (in particular, the internal pressure $p_c$) are considered uniform; in addition, the chamber temperature $T_c$ is considered constant in time. 

Nozzle throat erosion has been taken into account. A constant throat radius regression rate $r_t$ has been assumed in the present study, as usual during a preliminary performance estimation.\cite{Kamran2012} This allows for deriving analytic expressions for the consumed propellant mass (see next section).
Thus, the following relation can be exploited in order to determine the throat area value at a given time $t$ during motor operation:
\begin{equation}
     A_t (t) = A_{t,0} (1 + \overline{r}_t \eta)^2
 \label{eq:At}
\end{equation}
with $\overline{r}_t = \frac{E_t}{R_{t,0}}$ and $\eta = \frac{t}{t_{B,1}}$ a dimensionless regression rate and time, respectively, $E_t = R_{t,f} - R_{t,0}$ the total throat radius erosion, considered to be independent on the selected chamber pressure law, and $A_{t,0} = \pi R^2_{t,0}$ the initial throat area.

\subsection{Pressure Law Parameterization} 
\begin{figure} [h!]
    \centering
    \includegraphics[width = 0.6\textwidth]{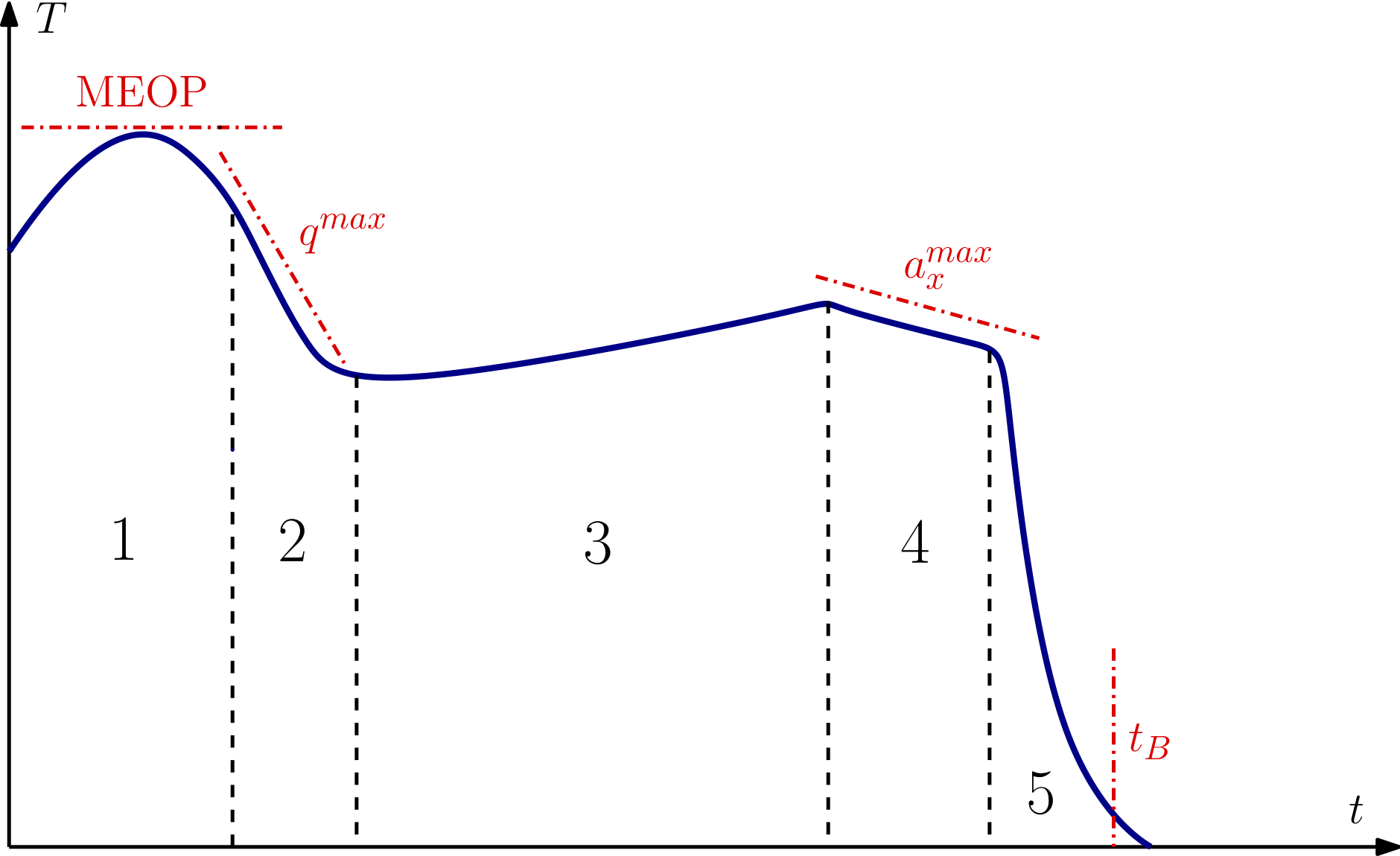}
    \caption{Typical first stage thrust profile}
    \label{fig:thrcurve}
\end{figure}

The thrust trend over-time for a first stage is mainly dictated by the several path constraints encountered along the atmospheric flight; complying with this constraints demands a particular shape for the thrust law $T(t)$, that, in turn, can be achieved by an adequate burning surface evolution $S_b(t)$.
Figure~\ref{fig:thrcurve} shows the main features of a typical first stage thrust profile determined by system requirements.

Five different regions can be clearly distinguished:
\begin{enumerate}
    \item 
    The maximum thrust value is attained soon, and limited by the maximum pressure $p_{meop}$ admitted within the combustion chamber, because of structural requirements, often called \textit{Maximum Expected Operating Pressure} (MEOP).
    \item Subsequently, the thrust decreases, in order to comply with the maximum dynamic pressure $q$ constraint.
    \item Once the altitude reached by the launcher becomes significant, the dynamic pressure value starts to drop; thus, an increase of the thrust is possible and beneficial.
    \item Near the end of the first stage flight, the vehicle mass is considerably reduced while the thrust is still quite high. A new decrease of the thrust is, therefore, the only way to bound the maximum value of the axial load $a_x$.
    \item The last part of the thrust law is the so called \textit{tail-off} phase; the rapid drop of pressure within the chamber occurs as a consequence of the sudden decrease of the grain burning surface. It may be advantageous to separate the engine and its associated dry mass before the propellant has been fully exhausted. 
\end{enumerate}

The optimal shape of the first stage thrust law is, here, attained through a parameterization of the chamber pressure law, which greatly resembles that of thrust. Acting on pressure has two clear advantages: (i) it allows, with respect to a direct thrust shaping, to account for the nozzle throat erosion and the effects of the expansion process; 
(ii) it maintains the computation cost low, with respect to a direct grain geometry optimization, and reduces the impact that necessary approximations for the SRM internal ballistic may have on the final result. 
The pressure profile will be modeled by patching together different curves, each represented by a simple time function and described by a small number of parameters. The motor burning time $t_{B,1}$ is considered as a free variable too. A parametric optimization process is exploited with the aim of determining the design variable values; the obtained profile is then re-scaled along pressure axis in order to comply with the assigned value of the total first stage propellant mass, being the Z40 loading capacity substantially fixed by the external case and central mandrel. In this way, only improvements due to the grain shape, rather than to propellant additional mass, are taken into account.

\subsubsection{PPLLLT Pressure Model.}
First of all, a dimensionless chamber pressure $\Pi = \frac{p}{p_{max}}$ and a dimensionless time $\eta = \frac{t}{t_{B,1}}$ are going to be used.

\begin{figure} [h]
    \centering
    \includegraphics[width = 0.7\textwidth]{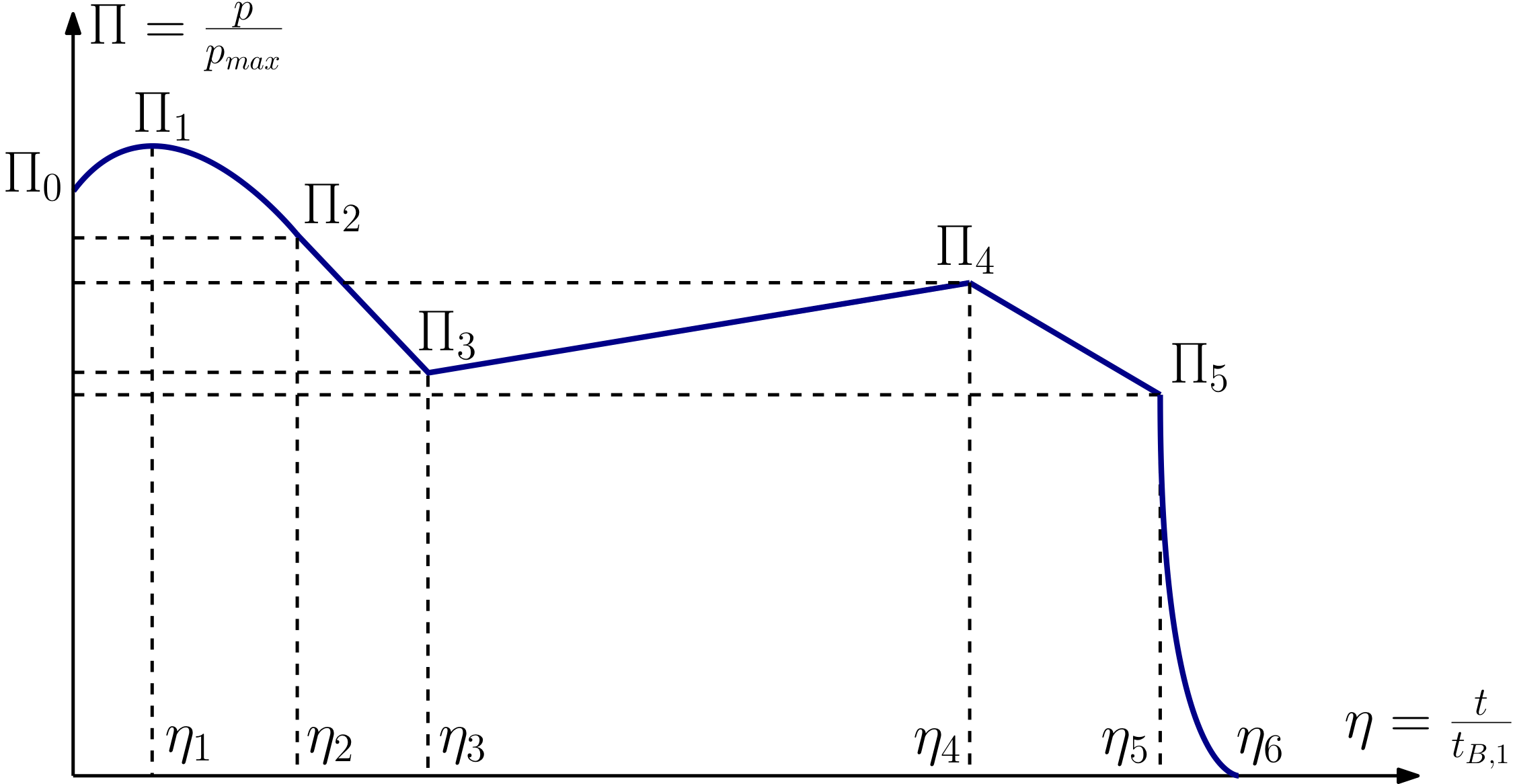}
    \caption{PPLLLT model for first stage chamber pressure.}
    \label{fig:PcPPLLLT}
\end{figure}

A pressure model that gives rise to a thrust law close to the typical one for a first stage (see Figure~\ref{fig:thrcurve}) is described by a Parabolic-Parabolic-Linear-Linear-Linear-Tailoff trend (here referred as \textit{PPLLLT model}).
The pressure time law is split into six consecutive legs: the first two are parabolic, simulating the 3-D burning of the star-shaped section of the finocyl grain, the next three are linear, each with its own slope, and the last one is, again, parabolic, approximating the tail-off phase, as shown in Figure~\ref{fig:PcPPLLLT}.

The dimensionless pressure trend is, therefore:
\begin{equation}
\begin{array}{lc}
      \Pi(\eta) = &
      \begin{cases}
      \frac{\Pi_0 - 1}{\eta^2_1}\,(\eta - \eta_1)^2 + 1 & \mbox{for}\,\,\, \eta \leq \eta_1\\
      \frac{\Pi_2 - 1}{(\eta_2 - \eta_1)^2}\,(\eta - \eta_1)^2 + 1 & \mbox{for}\,\,\, \eta_1 < \eta \leq \eta_2\\
       \Pi_{i-1} + \frac{\Pi_i - \Pi_{i-1}}{\eta_i - \eta_{i-1}}\, (\eta - \eta_{i-1}) & \mbox{for}\,\,\, \eta_{i-1} < \eta \leq \eta_i\,\,\, \forall i = 3,\dots,5\\
       a\, \eta^2 + b\,\eta + c & \mbox{for}\,\,\,  \eta_5 < \eta \leq 1
      \end{cases}
\end{array}
\label{eq:PcPPLLLT}
\end{equation}
being $\Pi_1 = \eta_6 = 1$ by their definition.

Pressure law during tail-off cannot be imposed arbitrarily, since it must reflect a precise behaviour which occurs while the last kilograms of propellant burn.
So, a quite accurate model for tail-off may be necessary to consider a realistic wait period to attain a relatively low thrust before stage separation occurs. Such wait period, in general, could have a significant effect on SRM performance.\cite{Panicker1998}
In particular, the initial steep descent of the chamber pressure during tail-off is due to a corresponding reduction of the burning surface, being the propellant next to finish. The final part of tail-off, instead, is often related to liner and thermal protections ablation phenomena\cite{Schiariti2015}; the effect of such residual thrust will be neglected in the present work. 
Experimental analysis conducted on a first-guess finocyl grain for Z40FS show that a parabolic law interpolates actual tail-off pressure with sufficient accuracy for the aim of the present study.
In order not to make the SRM behaviour during tail-off phase too much dependent on the obtained chamber pressure time law, it has been assumed that: (i) the same propellant mass $M_{p,to}$ is burnt during such a phase; (ii) the tail-off starts at the same fraction $\eta_{to}$ of the total burning time.

In order to guarantee that, during the search for optimum, the condition $\eta_{i - 1} \leq  \eta_i\,\, \forall i$ is always respected, each dimensionless time $\eta_i$ has been defined in terms of the previous time $\eta_{i - 1}$ and a fraction $\chi_{i}$ of the remaining time interval (until tail-off). On such basis, dimensionless times $\eta_i$ in Eq.~\eqref{eq:PcPPLLLT} are evaluated as:
\begin{equation}
\begin{array}{lc}
      \eta_i = &
     \begin{cases}
        \chi_1 & \mbox{for}\,\,\, i = 1 \\
        \frac{2(\Pi_2 - 1)\chi_3 + (\Pi_3 - \Pi_2)\chi_1}{2(\Pi_2 - 1)\chi_3 + (\Pi_3 - \Pi_2)} & \mbox{for}\,\,\, i = 2\\
        \eta_{i-1} + (\eta_5 - \eta_{i-1})\, \chi_i & \mbox{for}\,\,\, i = 3,\,4 \\
        \eta_{to} &  \mbox{for}\,\,\, i = 5\\
        1  & \mbox{for}\,\,\, i = 6
     \end{cases}
 \end{array}
\label{eq:etPPLLLT}    
\end{equation}
The dimensionless time $\eta_2$ has been evaluated so that it guarantees the derivative continuity at the corresponding time.

The total consumed propellant mass $M_p$, exploiting Eq.~\eqref{eq:MFR} and \eqref{eq:At}, is enforced as:
\begin{equation}
M_p =  \frac{t_{B,1}\, p_{max}\, A_{t,0}}{c^\ast}\, (m_{p,1} + m_{p,2} + m_{p,3} + m_{p,4} + m_{p,5}) + M_{p,to}
\label{eq:MpPPLLLT}
\end{equation}
where:
\begin{equation}
\begin{array}{lc}
  m_{p,i} = \int_{\eta_{i - 1}}^{\eta_i} \Pi(\eta)\,  (1 + \overline{r}_t\, \eta)^2 d\eta  & \mbox{for}\,\,\, i = 1,\dots,5 \\ 
\end{array}
\label{eq:mpi}    
\end{equation}
Analytic expressions for $m_{p,i},\,\forall i=1,\dots,5$, as a function of the pressure law design variables and $\overline{r}_t$ can be easily obtained. 
At this point, the maximum value for the operating pressure $p_{max}$ can be determined:
\begin{equation}
 p_{max} = \frac{c^\ast}{t_{B,1}\, A_{t,0}}\, \frac{M_p\, -\, M_{p,to}}{m_{p,1} + m_{p,2} + m_{p,3} + m_{p,4} + m_{p,5}}
\label{eq:pmeopPPLLLT}
\end{equation}
It should be remarked that the constraint on the MEOP must be respected:
\begin{equation}
p_{max} \leq p_{meop}
\label{eq:ConMEOP}    
\end{equation}
Because of the assumed hypotheses, the tail-off parabolic law can be determined by imposing the following conditions:
\begin{equation}
    \begin{split}
          \Pi(\eta_5) &= \Pi_5\\
          \Pi(1) &= 0\\
          \int_{\eta_5}^{1} \Pi(\eta) d\eta &= m_{p,6}
    \end{split}
\label{eq:TO1}
\end{equation}
where $m_{p,6}$ is equal to: 
\begin{equation}
        m_{p,6} = \frac{M_{p,to} c^\ast}{p_{max}\, t_{B,1}\, A_{t,0}\, (1 + \overline{r}_t\, \eta_5)^2}
\label{eq:mpV}
\end{equation}
In particular, for the sake of simplicity, it has been assumed that, because of the fast drop-off of mass flow rate during tail-off, the nozzle throat erosion is negligible in such phase (\emph{i.e.}, $r_t \equiv 0\,\,\, \mbox{for}\,\,\, \eta \geq \eta_5$).
Consequently, the following expressions for the the tail-off parabolic law coefficients result:
\begin{equation}
\begin{split}
         a &= \frac{3}{(1 - \eta_5)^2}\left(\Pi_5 - 2\frac{m_{p,6}}{(1 - \eta_5)}\right) \\
         b &= -\frac{\Pi_5}{1 - \eta_5} - a(1 + \eta_5)\\
         c &= - a - b
\end{split}
\end{equation}
Not all the values for $\Pi_5$ are, however, admissible. Indeed, so that the parabolic trend correctly approximates the real pressure trend during tail-off, the following constraints, acting on $\Pi_5$, must be set:
\begin{equation}
    \begin{split}
         a \geq 0 &\implies \Pi_5 \geq \frac{2 m_{p,6}}{1 - \eta_5} \\
         \eta_{-} = 1 &\implies \Pi_5 \leq \frac{3 m_{p,6}}{1 - \eta_5}
    \end{split}
\end{equation}
The former guarantees a convex tail-off parabola, the latter places the first zero $\eta_{-}$ of the parabola in $\eta_{6} = 1$.

The design variables referred to the first stage chamber pressure with their respective boundaries, are listed in Table~\ref{tab:MotorLaws}.
Particular attention must be paid in the selection of the boundaries of each variable. In fact, such boundaries must be:
(i) sufficiently wide, so that the search space extension is not limited too much, eventually ruling out some potential optimal solutions; (ii) sufficiently close, to avoid that the algorithm finds putative optimal solutions that, in practise, cannot result in a realistic initial geometry of the grain.

\begin{table}[h!]
    \centering
   \caption{PPLLLT model for Z40FS chamber pressure.}
    \label{tab:MotorLaws}
    \begin{tabular}{c|c|c|c}
      \hline 
       \textbf{Model} & \textbf{Time} & \textbf{Optimization} & \textbf{Variables}\\
        & \textbf{law} & \textbf{variables} & \textbf{boundaries}\\
        \hline
       \textbf{PPLLLT} & \textbf{parabolic} & $\Pi_{0}$ & $[0.85,\,0.95]$\\
       & \textbf{parabolic} & $\chi_1$ & $[0.05,\,0.15]$\\
       & \textbf{linear} & $\Pi_2$ & $[0.8,\,0.95]$\\
       & \textbf{linear} & $\chi_3$ & $[0.1,\,0.3]$ \\
       & \textbf{linear} & $\Pi_3$ & $[0.6,\,0.75]$\\
       & \textbf{tail-off} & $\chi_4$ & $[0.4,\,0.85]$\\
       &  & $\Pi_4$ & $[0.6,\,0.75]$\\
       & & $\Pi_5$ & $[0.6,\,0.75]$\\
      \hline 
    \end{tabular}
\end{table}

\section{Optimization}

\subsection{Optimization Problem}
Solving an optimization problem means to find a set of design variables $\mathbf{x} = [x^{(1)}, \dots , x^{(N_D)}]$ which maximize an objective function $f(\mathbf{x})$.
The feasible solution space, \emph{i.e.}, the space within which it is admitted to look for the optimal design variables, is often limited by:\newline
(i) the presence of an eligibility interval for each variable $x^{(i)}$ of the solution vector $\mathbf{x}$:
\begin{equation}
    x^{(i)}_L \leq x^{(i)} \leq x^{(i)}_U\,\,\,\,\,\,\,\, \mbox{for}\,\,\, i = 1,\dots,N_D
\label{eq:LBUB}
\end{equation}
(ii) the existence of inequality constraints of the form:
\begin{equation}
    \Psi^{(j)}(\mathbf{x}) \leq 0\,\,\,\,\,\,\,\, \mbox{for}\,\,\, j = 1,\dots,k
\label{eq:ineqCon}
\end{equation}
and, in this case, one is dealing with a \textit{constrained optimization problem}. 

In the \textit{launcher trajectory optimization problem} that is going to be faced, the objective is to maximize the transported payload: $f(\mathbf{x}) = M_u$. The problem can be formulated as follows:
\begin{equation}
\begin{cases}
     \max\limits_{\mathbf{x}_L \leq \mathbf{x} \leq \mathbf{x}_U}  M_{u} \\
     \mbox{s.t.}\,\,\,  \Psi^{(j)}(\mathbf{x}) \leq 0 & \mbox{for}\,\,\, j = 1,\dots,k
\end{cases}
\label{eq:optProbasc}   
\end{equation}
where: 
\begin{equation}
     \mathbf{x} = M_u \cup \mathbf{x}_{flight} \cup \mathbf{x}_{guidance} \cup \mathbf{x}_{motor} 
\label{eq:xvec}   
\end{equation}
In particular:
\begin{align}
     \mathbf{x}_{flight} &= [t_{1,2}, t_{1,3}, t_{B,1}, t_{C2}, t_{C4}] \label{eq:xflight}\\
     \mathbf{x}_{guidance} &= [ \theta_{T}(t_2), \psi_{T}(t_2), \theta_{O}(t_7), \theta_{O}(t_8), \xi, \psi_{T}(t_7), \theta_{O}(t_9), \psi_{O}(t_9)] \label{eq:xguidance}  \\
     \mathbf{x}_{motor} &= [\Pi_0,\, \chi_1,\, \Pi_2,\, \chi_3,\, \Pi_3,\, \chi_4,\, \Pi_4,\, \Pi_5] \label{eq:xmotor}
\end{align}
Lower and upper bounds for flight/guidance variables and motor variables has been defined in Table~\ref{tab:GuidLaws} and Table~\ref{tab:MotorLaws}, respectively.
Concerning the constraints, these are: \newline \newline
\textit{Path Constraints}:
\begin{align}
\Psi^{(1)} &= q - q^{max} & \mbox{dynamic pressure}\\
\Psi^{(2)} &= (q \ast \alpha) - (q \ast \alpha)^{max} & \mbox{bending load}\\
\Psi^{(3)} &= \Dot{Q}_W - \Dot{Q}_W^{max} & \mbox{heat flux at fairing jettisoning}\\
\Psi^{(4)} &= a_x - a^{max}_{x} & \mbox{axial acceleration}\\
\Psi^{(5)} &= \dot{\theta}_{PO} - \dot{\theta}^{max}_{PO} & \mbox{\textit{pitch-over} rate}
\end{align}
\textit{Terminal Constraints}:
\begin{align}
&
\begin{cases}
\Psi^{(6)} = r_a - (\tilde{r} + \delta r) \\
\Psi^{(7)} = - [r_p - (\tilde{r} + \delta r)]
\end{cases} 
& \mbox{assigned orbital radius}\\
&\,\,\,\,\,\, \Psi^{(8)} = |i - \tilde{i}| - \delta i & \mbox{assigned orbital inclination}
\end{align}
\textit{Pressure Law Constraints}:
\begin{align}
\Psi^{(9)} &= p_{max} - p_{meop} & \mbox{MEOP}\\
\Psi^{(10)} &= - (\Pi_5 - \frac{2 m_{p,6}}{1 - \eta_5}) & \mbox{convex tail-off parabola}\\
\Psi^{(11)} &= \Pi_5 - \frac{3 m_{p,6}}{1 - \eta_5} & \mbox{parabola first zero at $t_{B,1}$} 
\end{align}

\subsection{Global Optimization Algorithm}
An optimization algorithm based on \textit{Differential Evolution} (DE), that proved to be effective in other space trajectory optimization problems\cite{Federici2018}, has been selected in the present application. DE is a population-based stochastic meta-heuristic algorithm, firstly introduced by R. Storn and K. Price in 1997,\cite{Storn1997} featuring simple and efficient heuristics for global optimization problems defined over a continuous space. Being inspired by evolution of species, it exploits the operations of \textit{Cross-over}, \textit{Mutation} and \textit{Selection} to generate new candidate solutions. Because of its good performance on several benchmarks and real-world optimization problems, many researchers all over the world have directed their efforts to further improve the effectiveness of the original algorithm, by devising many variants, collected in Reference~\citenum{Das2011}. 

In the present implementation four different mutation strategies, among those originally proposed\cite{Storn1997}, in conjunction with a binomial-type crossover, have been adopted. By using the current nomenclature\cite{Storn1997,Das2011}, the DE schemes here exploited are:
\begin{enumerate}
    \item \textit{DE/rand/1/bin}
    \item \textit{DE/best/1/bin}
    \item \textit{DE/target-to-best/1/bin}
    \item \textit{DE/best/2/bin}
\end{enumerate}
In particular, strategies based on mutation of the best individual (strategies 2 and 4) typically show a faster converge rate toward an (often local) minimum, whereas strategies based on randomly chosen individuals (strategies 1 and 3) better explore the
whole search space.
In order to avoid the manual tuning of DE control parameters, a self-adaptation scheme\cite{Brest2007} has been implemented: the values of the scale factor $F$ and the crossover probability $C_r$ are encoded into the individuals that undergo the optimization procedure.

\begin{figure} [h!]
    \centering
    \includegraphics[width = 0.7\textwidth]{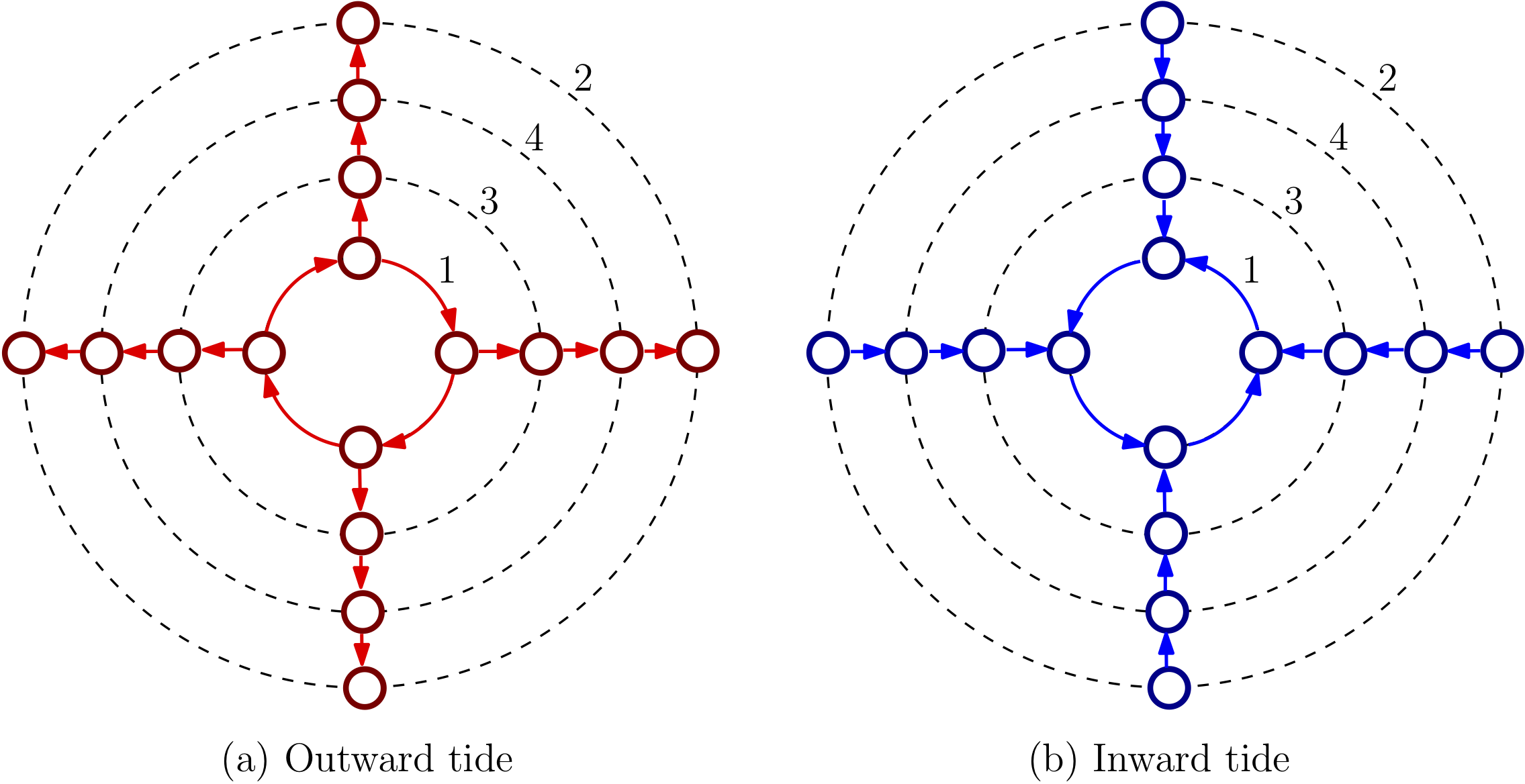}
    \caption{Migration scheme: forward (a) and backward (b), for the 16-island case.}
    \label{fig:Islands}
\end{figure}

In order to achieve a good balance between the search space exploration and a faster convergence to a good solution, the proposed algorithm involves the creation of different populations, or ``tribes'', each located on an ``island'' of an archipelago, arranged in a radial configuration, as Figure~\ref{fig:Islands} shows. Each tribe evolves independently from the others and features one specific mutation strategy among the four proposed variants (as reported in Figure~\ref{fig:Islands}) until a migration is performed, during which each tribe passes its best agents to the ``following one'' (according to the direction of migration), in which these agents replace the same number of worst agents. Migration tides alternate at each event. The proposed island-model allows an easy parallelization\cite{Izzo2009} on architectures with $4\times$ cores. The termination criterion of the optimization process is simply based on the generation number. A partial-restart mechanism, named ``Epidemic'', has been adopted to maintain diversity between individuals in each population and avoid a premature convergence on a local optimum.

\subsubsection{Constraint Handling.}
Evolutionary algorithms are, in their ``standard'' formulation, unconstrained optimization methods; this fact makes it difficult to understand how to efficiently incorporate constraints of any kind into the objective function during the search.
The $\varepsilon$DE approach proposed by T. Takahama and S. Sakai\cite{Takahama2010}, obtained by applying the $\varepsilon$ constrained method to differential evolution, is here adopted.
The idea is to decide the selection process through a set of simple rules:
\begin{enumerate}
    \item between feasible individuals, the one with the maximum value of the objective function wins;
    \item between a feasible and an infeasible individual, the feasible one is selected;
    \item between two infeasible individuals, the one with the lowest value of the maximum constraint violation:
    \begin{equation}
        \Psi^{max} = \max\limits_{j \in [1, k]} \Psi^{(j)}(\mathbf{x})
    \end{equation}
    survives in the next generation.
\end{enumerate}
Without further refinements, this simple approach pushes the evolution process to focus first on the constraint satisfaction and, then, to randomly sample points in the identified feasible region.\cite{Coello1999}
Therefore, a constraint satisfaction tolerance $\varepsilon$ for $\Psi^{max}$ is introduced and decreased during generations $G$ according to the rule:
\begin{equation}
\begin{array}{cc}
     \varepsilon^{G} = &  
\begin{cases}
 \varepsilon^{0} & \mbox{for}\,\,\, G \leq N^{0}\\
 \varepsilon^{0} \left[\frac{\varepsilon^{\infty}}{\varepsilon^{0}}\right]^{\frac{G - N^{0}}{N^{\infty} - N^{0}}} & \mbox{for}\,\,\, N^{0} < G < N^{\infty}\\
 \varepsilon^{\infty} & \mbox{for}\,\,\, G \geq N^{\infty}
\end{cases}
\end{array}
\label{eq:eps}    
\end{equation}
with $\varepsilon^{0}, \varepsilon^{\infty}$ the initial and final value of the tolerance, respectively.
In this way, ``bad'' moves (i.e. toward infeasible solutions) are more likely to be allowed at the beginning, when the tolerance $\varepsilon$ is high and the entire search space must be explored in order to identify promising regions; as the generation number increases, such moves became forbidden, and the search concentrates in the feasible part of the identified ``optimal'' region.
The following values for the above parameters have been used in the present application: $\varepsilon^{0} = 0.05$, $\varepsilon^{\infty} = 1\times 10^{-8}$, $N^{0} = \frac{1}{6}N_G$, $N^{\infty} = N_G$.

\section{Results}
The reported solution has been obtained through an 8-island optimization engine, with $N_p = 50$ agents per tribe and a maximum number of generations equal to $N_G = 5000$. 
In order to increase the confidence on the attained result, the optimization is repeated several times, and the best found solution is taken as putative optimum.
A circular sun-synchronous orbit with an altitude of $500$ km above Earth surface and an inclination of $98\degree$ has been considered as target mission.

The results of the optimization process, in term of optimal values for the design variables, are reported in Table~\ref{tab:SSO_PPLLLT1}. Table~\ref{tab:SSO_PPLLLT2}, instead, lists: the initial and final inertial (ECI) velocity of the vehicle $v_i$ and $v_f$, the velocity increment provided by the propulsion system $\Delta v_{prop}$, and the cumulative amount of velocity losses (gravitational $\Delta v_{grav}$, aerodynamic $\Delta v_{aero}$ and misalignment $\Delta v_{mis}$) along the obtained ascent trajectory.

\begin{table}[h!]
    \centering
    \caption{Optimum solution vector $\mathbf{x}$.}
    \begin{tabular}{c|S|c|S|c|S}
      \hline
      \multicolumn{6}{c}{$ M_u$ $= 409\,{kg}$}\\
      \hline
       \multicolumn{2}{c|}{$\mathbf{x}_{flight}$}  & \multicolumn{2}{c|}{$\mathbf{x}_{guidance}$} & \multicolumn{2}{c}{$\mathbf{x}_{motor}$} \\
       \hline
        $ t_{1,2}$ $[s]$ & 9.14 & $ \theta_T(t_2)$ [\degree] & 71.72 & $ \Pi_0 $ & 0.850 \\
       $ t_{1,3}$ $[s]$ & 4.98 & $ \psi_T(t_2)$ [\degree] & 190.14 & $ \chi_1 $ & 0.142\\
       $ t_{B,1}$ $[s]$ & 95.47  & $ \theta_O(t_7)$ [\degree] &  6.34 & $ \Pi_2 $ & 0.923\\
       $ t_{C2}$ $[s]$ & 89.55 & $ \theta_O(t_8)$ [\degree] & 0.40 & $ \chi_3 $ & 0.102\\
       $ t_{C4}$ $[s]$ & 2451.25 & $ \xi$ & -0.21 & $ \Pi_3 $ & 0.721\\
       $ - $ & & $ \psi_O(t_7)$ [\degree] & -4.48 & $ \chi_4$ & 0.834\\
       $ - $ & & $ \theta_O(t_9)$ [\degree]& -1.09 & $ \Pi_4 $ & 0.749\\
       $ - $ & & $ \psi_O(t_9)$ [\degree] & 0.81 & $ \Pi_5 $  & 0.652\\
    \hline
    \end{tabular}
    \label{tab:SSO_PPLLLT1}
\end{table}

\begin{table}[h!]
    \centering
    \caption{Initial and final inertial velocity $\bm{v_i}$ and $\bm{v_f}$, and velocity variations $\bm{\Delta v}$ along the trajectory.}
    \begin{tabular}{c|c|c|c|c|c}
      \hline
       $v_i$ & $v_f$ & $\Delta v_{prop}$ & $\Delta v_{grav}$ & $\Delta v_{aero}$ & $\Delta v_{mis}$\\
       $[km/s]$ & $[km/s]$ & $[km/s]$ & $[km/s]$ & $[km/s]$ & $[km/s]$ \\
       \hline
       0.37 & 7.61 & 9.46 & 1.39 & 0.14 & 0.68\\
     \hline
        \end{tabular}
    \label{tab:SSO_PPLLLT2}
\end{table}

\begin{figure} [h!]
    \centering
    \includegraphics[width = 0.8\textwidth]{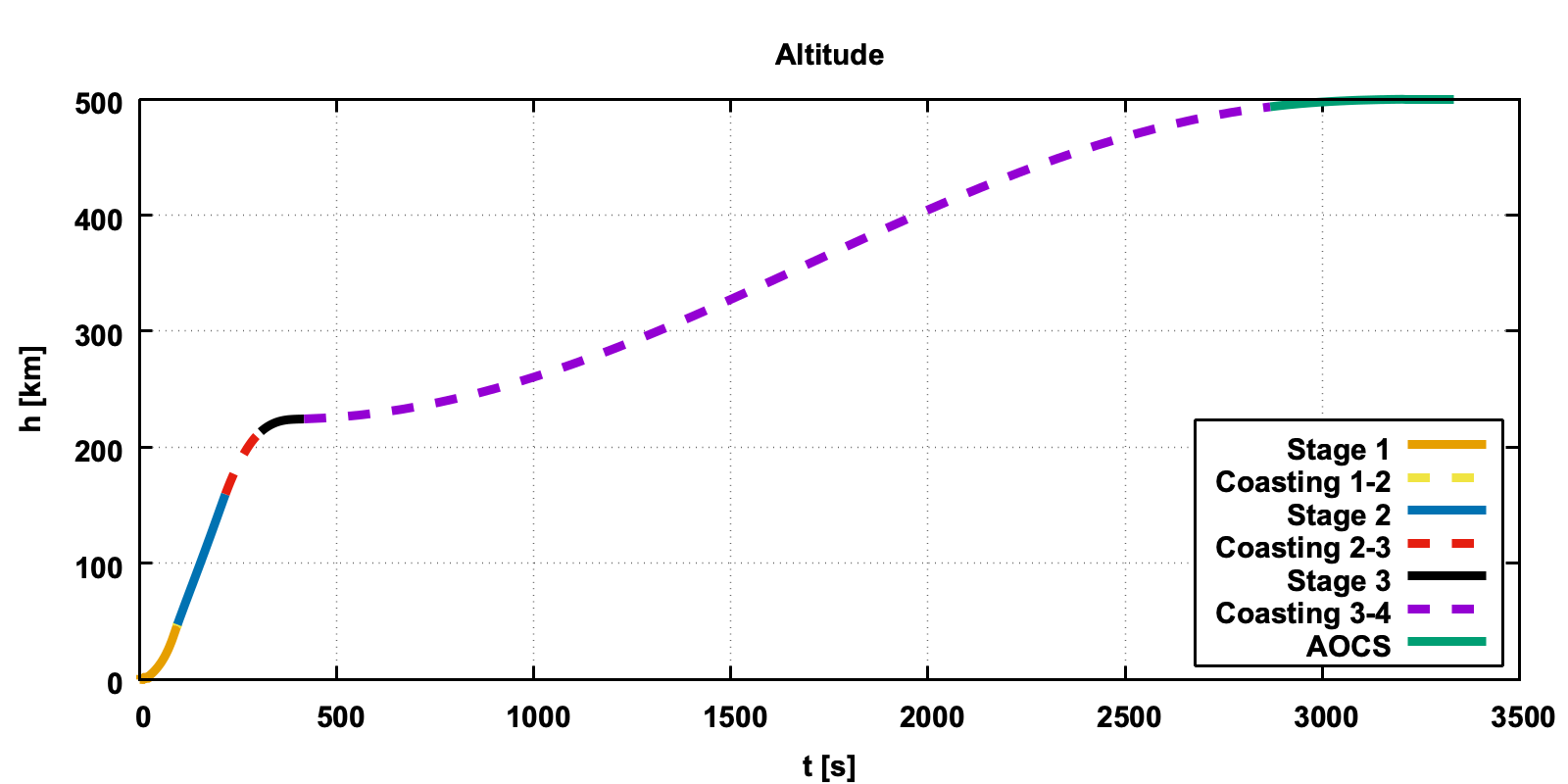}
    \caption{Altitude vs time}
    \label{fig:AltPPLLLT}
\end{figure}

Figure~\ref{fig:AltPPLLLT} shows the altitude profile along the trajectory. A direct ascent trajectory is flown until the third stage ignition, which occurs just outside Earth atmosphere. At this point, the velocity gained during ZX operation enters the vehicle at the perigee of an Hohmann-like (because of the finite thrust) ellipse towards the target orbit. A \textit{coasting 3-4} of half an orbit period lets the rocket reach the apoapsis, where the final orbital circularization is carried out by means of the ZX AOCS.
This is the typical mission profile for a medium-altitude final orbit; a complete direct ascent would, indeed, entail high gravitational losses, because of the great value of the flight angle between inertial velocity and the horizon.
It should be noticed that, actually, the considered AOCS is capable of a unique shoot and ZX must cover this deficiency by realizing the first $\Delta v$ of the Homhann-like transfer; on the contrary, during a typical mission of the \textit{Vega} launcher, the Hohmann maneuver is completely carried out by AVUM, able to perform multiple ignitions.




\begin{figure} [h!]
    \centering
    \includegraphics[width = 0.8\textwidth]{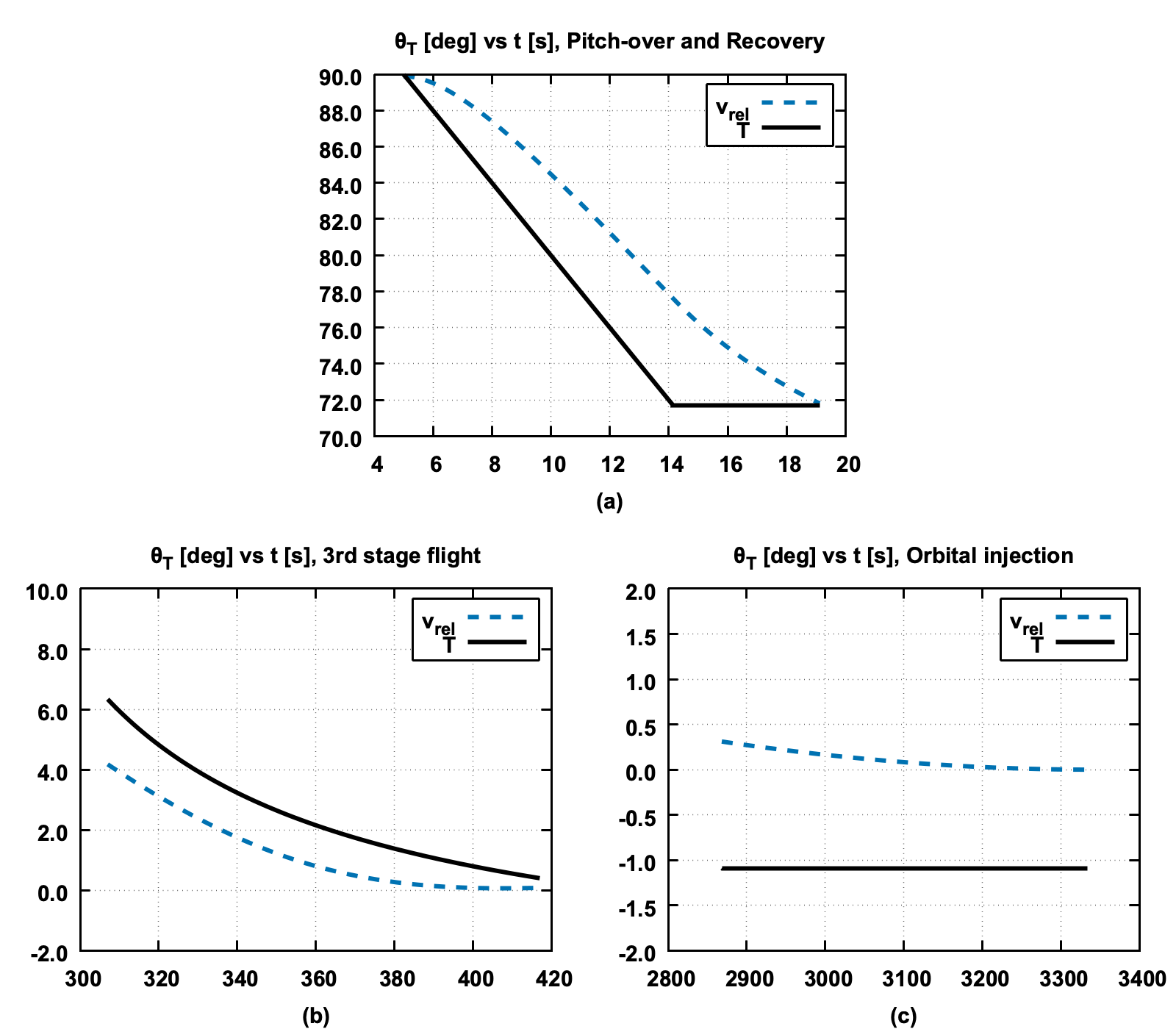}
    \caption{Flight path angle $\bm{\theta_T}$ vs time, during: (a) \textit{pitch-over} and \textit{recovery}, (b) \textit{3rd stage flight}, (c) \textit{orbital injection}.}
    \label{fig:angles_PPLLLT}
\end{figure}

Figure~\ref{fig:angles_PPLLLT} shows the topocentric flight path angle during \textit{pitch-over} and \textit{recovery} (a), \textit{3rd stage flight} (b) and final \textit{orbital injection} (c). 
The linear pitch law and the constant attitude used during \textit{pitch-over} and \textit{recovery}, respectively, are clearly visible in image (a); a perfect alignment between thrust and relative velocity occurs just before entering the first gravity turn maneuver, as desired. The bi-linear tangent law performed during ZX flight is reported in image (b), whereas the constant orbital flight path angle assumed during AOCS operation is shown in image (c). 
From Figure~\ref{fig:angles_PPLLLT} one can appreciate that the topocentric flight path angle referred to launcher relative (and, then, inertial) velocity reaches a null value both at the end of the \textit{3rd stage flight} and at the end of the \textit{coasting 3-4}: this is a clear confirmation that a Hohmann-like transfer brings the vehicle to its final orbit, after ZX shutdown. 

\begin{figure} [h!]
    \centering
    \includegraphics[width = 0.7\textwidth]{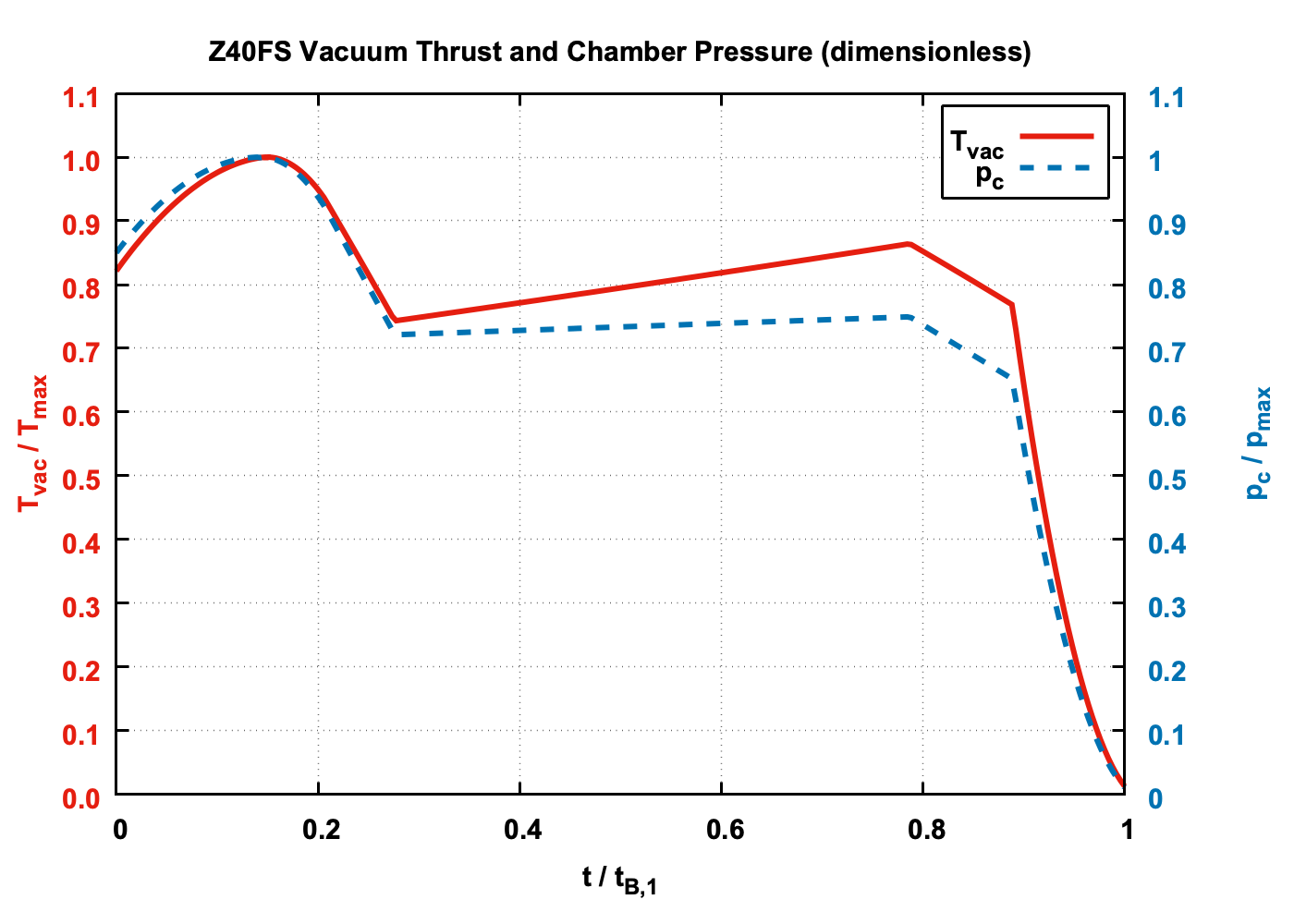}
    \caption{First-stage thrust (red curve) and chamber pressure (light blue curve) vs time.}
    \label{fig:ThrPc_PPLLLT}
\end{figure}


Thrust and pressure laws for the first stage are depicted in Figure~\ref{fig:ThrPc_PPLLLT}. The depicted quantities, as well as the time axis, have been re-scaled so that they belong to interval $[0, 1]$. The great similarity with the typical first stage thrust (see Figure~\ref{fig:thrcurve}) immediately catches the eye.

\begin{figure} [h!]
    \centering
    \includegraphics[width = 0.8\textwidth]{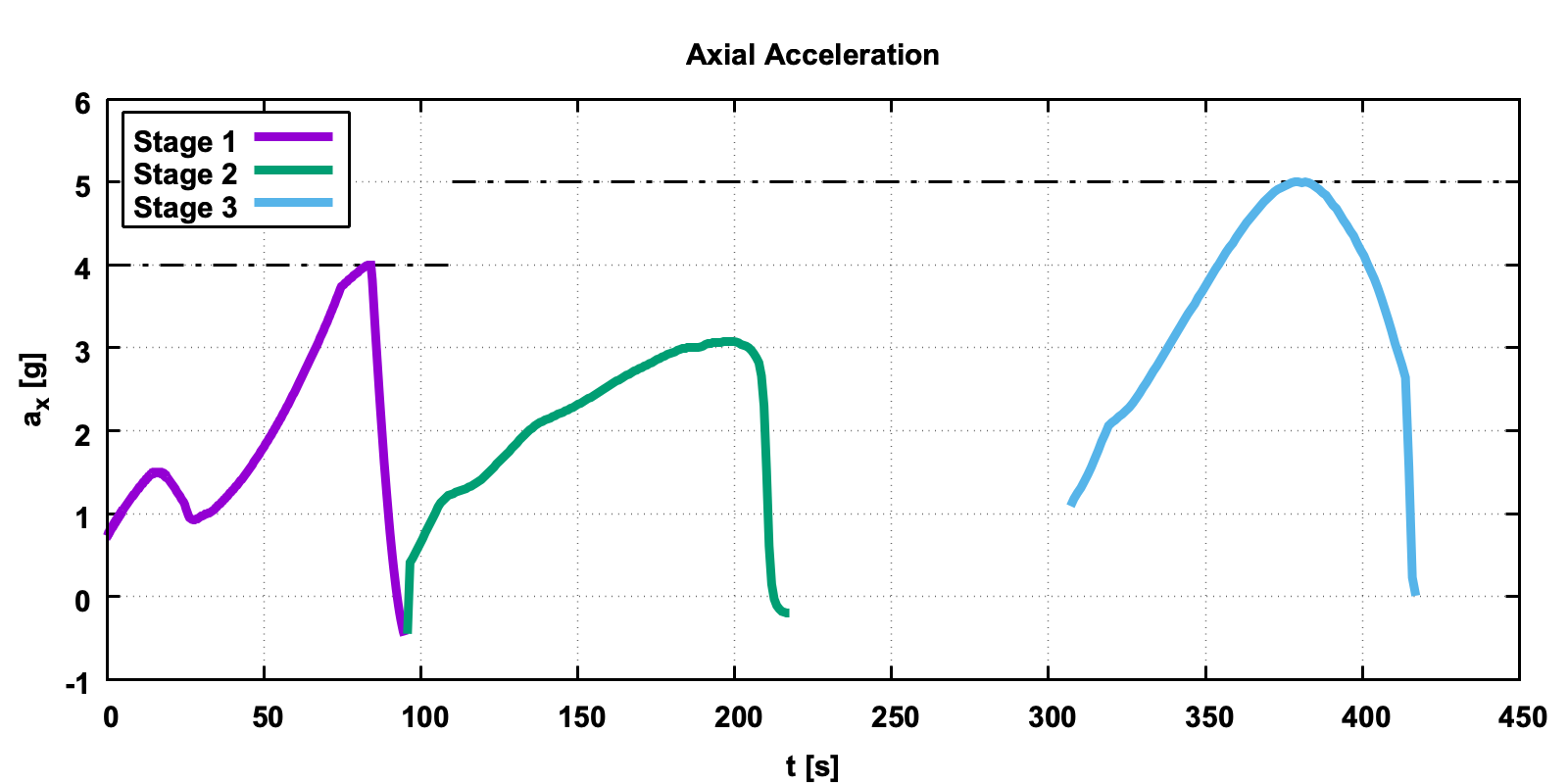}
    \caption{Axial acceleration vs time.}
    \label{fig:AccPPLLLT}
\end{figure}

Through the initial parabolic trend the thrust reaches its maximum value very early along the ascent (\emph{i.e.}, during \textit{pitch-over}); the great acceleration lets the vehicle acquire, very soon, a relatively high velocity, reducing the amount of gravitational losses (it is, indeed, evident in Table~\ref{tab:SSO_PPLLLT2} how important gravitational losses are among the velocity losses). 
\begin{figure} [h!]
    \centering
    \includegraphics[width = 0.8\textwidth]{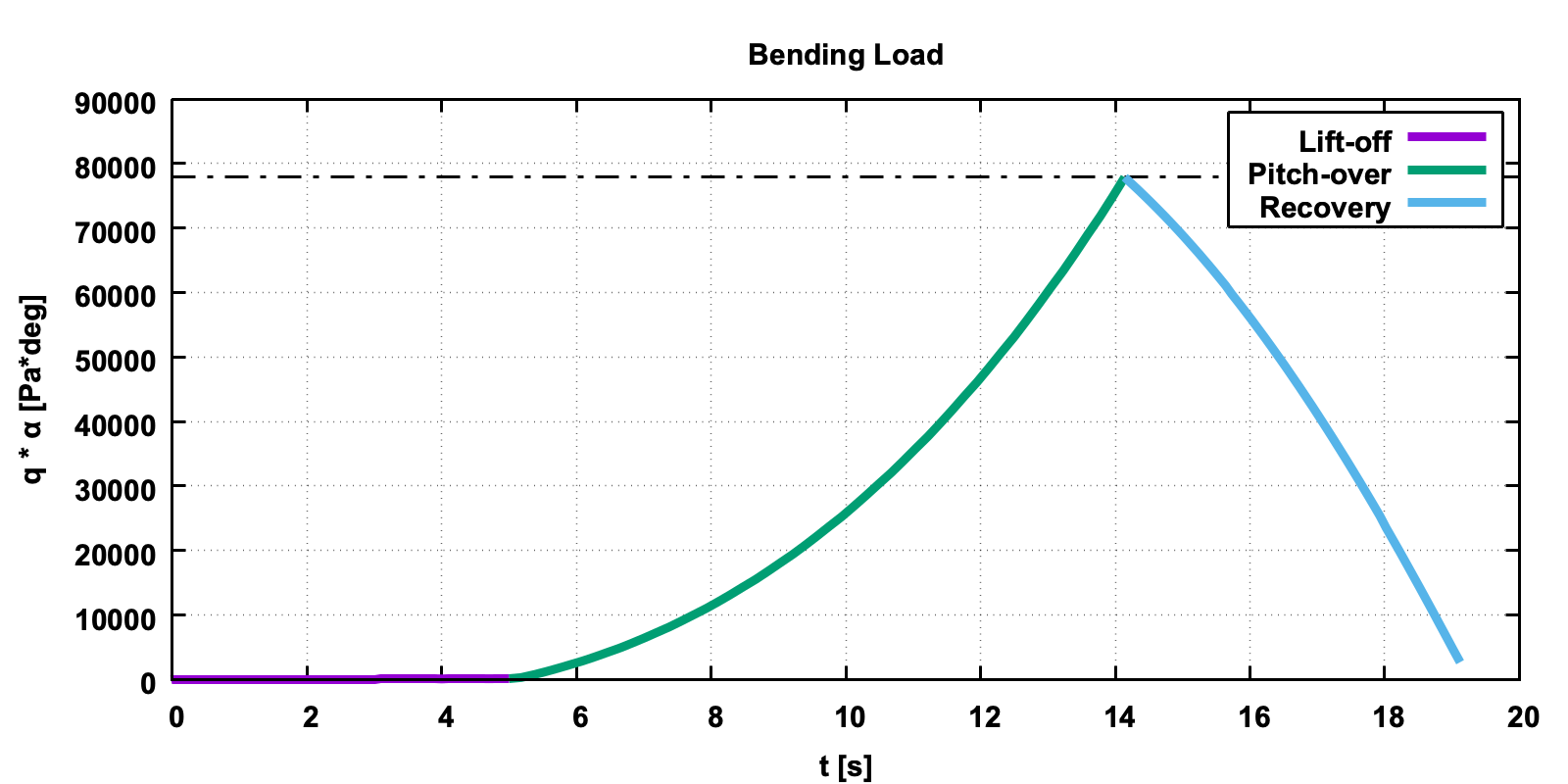}
    \caption{Bending load vs time.}
    \label{fig:qalphaPPLLLT}
\end{figure}
The condition of maximum bending load is reached at this point of the flight; as Figure~\ref{fig:ThrPc_PPLLLT} and Figure~\ref{fig:qalphaPPLLLT} clearly display, the thrust value is still quite high while $(q \ast \alpha)$ touches its limit value; it is, indeed, convenient to limit the magnitude of the lift force by decreasing the angle of incidence $\alpha$ rather then by reducing the rocket acceleration, fundamental in these initial phases of the ascent to maintain gravitational losses low, and, consequently, gain in payload mass. The non-null bending load $(q \ast \alpha)$, together with the high misalignment and aerodynamic losses, pushes the \textit{pitch-over} rate toward its maximum admitted value: in this way, the desired \textit{kick-angle} is attained with the minimum \textit{pitch-over} time-length.
\begin{figure} [h!]
    \centering
    \includegraphics[width = 0.8\textwidth]{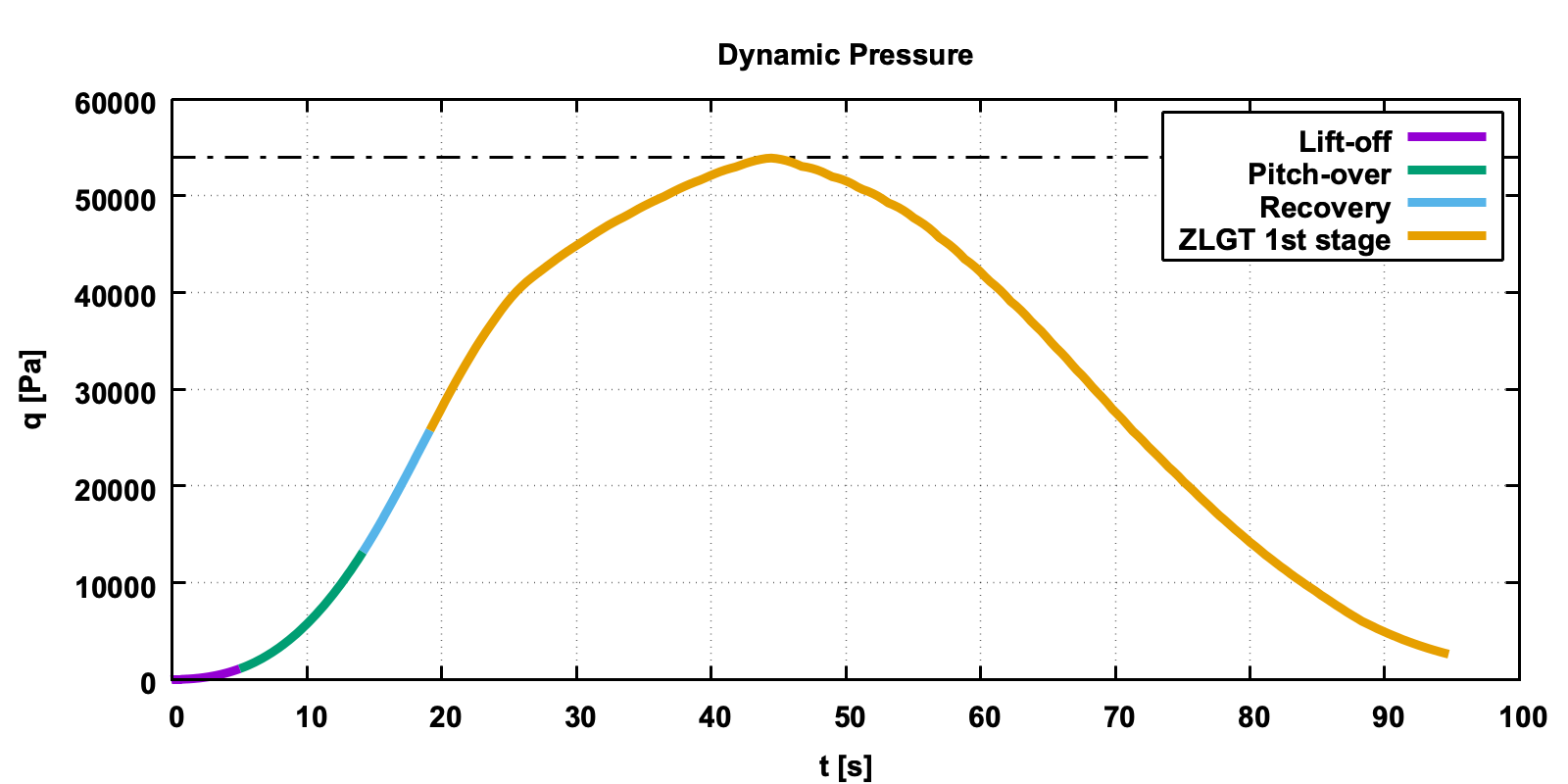}
    \caption{Dynamic pressure vs time.}
    \label{fig:pdynPPLLLT}
\end{figure}
The parabolic region is followed by a linear decreasing law, dictated by the maximum value admitted for dynamic pressure. It is easy to note from Figure~\ref{fig:pdynPPLLLT} the typical increasing-decreasing trend of the dynamic pressure, as a direct result of the exponential reduction in atmospheric density with altitude and the increasing velocity of the launcher. The so-called \textit{maximum q} condition is encountered during the \textit{ZLGT 1st stage} phase and approximately occurs at the beginning of the increasing linear region of pressure law; once such condition is left behind, the thrust level is free to grow again.
As the propellant mass is approaching its end, the axial acceleration climbs so rapidly that it is necessary to stop the thrust increment in order not to violate the corresponding constraint (see Figure~\ref{fig:AccPPLLLT}).
A final parabolic tail-off drives the thrust toward its minimum; Z40 is detached from the rest of the vehicle at the time the pressure reaches an imposed minimum threshold $p_{thre} = 1\,bar$.
It is possible to note that dynamic pressure, bending load and axial acceleration reach their maximum admitted value, depicted in graphs as a dashed straight line: this is a clear evidence of the optimality of the solution found so far and of the good performance of the optimization algorithm in terms of constraint handling. 
Indeed, it is foreseeable that the best performance in terms of payload mass should be obtained when the launcher works at its maximum structural capacity, thus exploiting at best its propulsive capability along a trajectory that limits losses.
The maximum heat flux after fairing jettisoning, instead, is well below its threshold value; this suggests that fairing ejection could be performed during the second stage flight, maybe allowing for a further improvement of the payload.





\begin{figure} [h!]
    \centering
    \includegraphics[width = \textwidth]{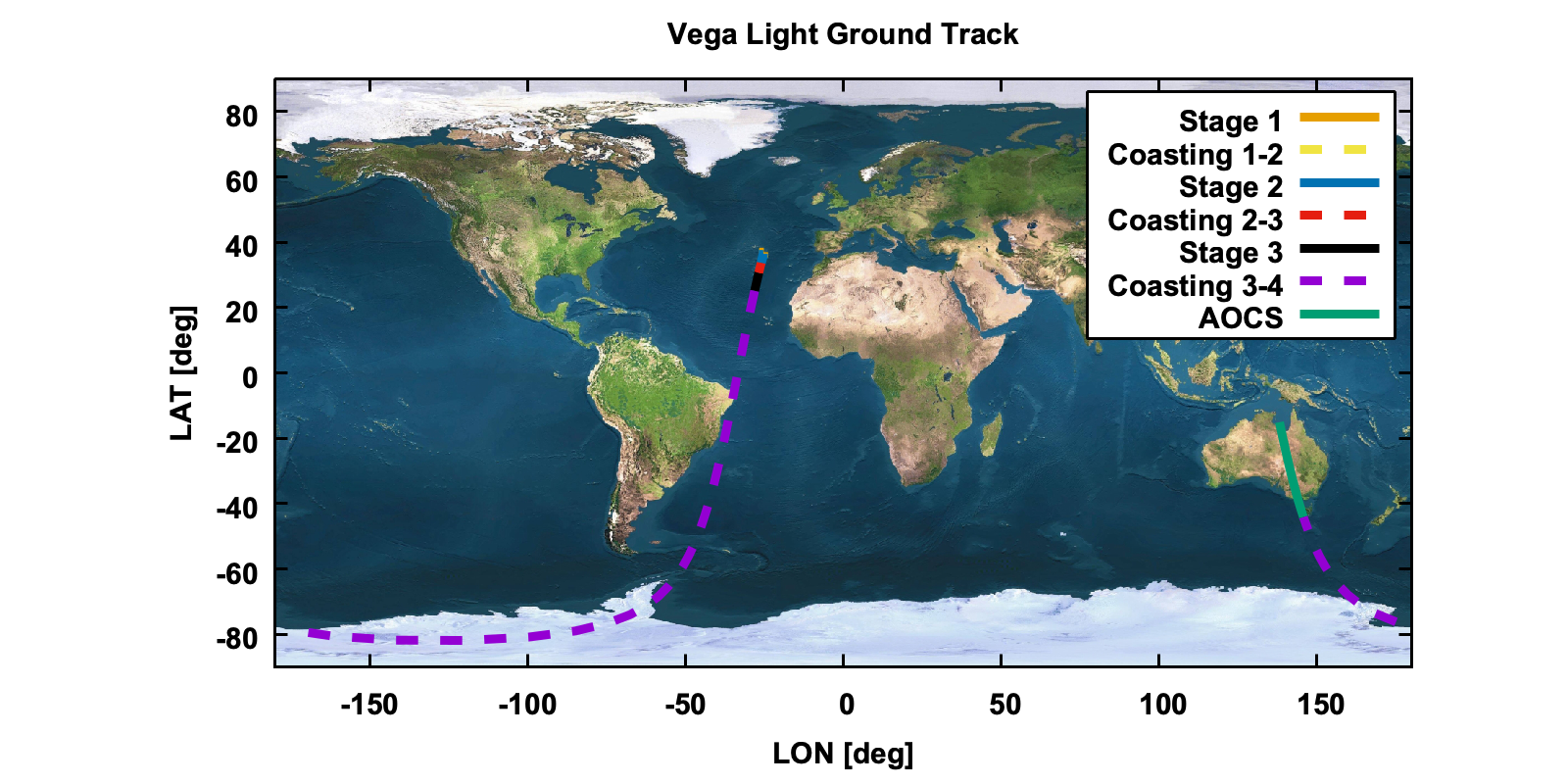}
    \caption{Launcher ground track.}
    \label{fig:GT_PPLLLT}
\end{figure}

Even though, in the present analysis, the geographical constraints concerning the fall-down points of the first two stages have not been analyzed, it can be noticed that a south-ward launch from Azores' base allows for an initial flight over the Atlantic Ocean, where no land is encountered, as shown in Figure~\ref{fig:GT_PPLLLT}.

\section{Conclusion}
This paper proposes an optimization methodology that offers the great advantage, over ``classical'' ascent trajectory optimization techniques, to account for the close and mutual relationship that exists between SRM thrust law and flight mechanics. 
The methodology is quite general and can be implemented, with the appropriate changes related to the launcher configuration, every time the determination of the optimal thrust of first-stage main motors or strap-on boosters is required, in addition to the mandatory optimization of the rocket trajectory, without the need to perform a complex analysis of the 3-dimensional SRM internal ballistic and of the real regression rate of the propellant grain surface. As a post-process step, further studies have to be performed with the aim of determining the real grain geometry that produces the previously-obtained chamber pressure history. The analysis can be extended also to upper stage SRM grain design, although the strict path constraints which characterize the atmospheric flight of the launcher first stage have a great influence in modelling the shape of its thrust law. 

The chamber pressure time law has been modeled using few analytic functions; this allows for an easy enforcement of the total propellant mass to be loaded into the motor, avoiding the use (i) of additional constraints, which would further limit the search space, and (ii) of quadrature formulae for the approximation of mass flow rate integrals, which worsen the computational load and might bring to dissimilarity in the order of tens of kilograms of propellant with respect to the imposed mass, because of their limited accuracy.
Despite the exposed optimization methodology is still at its beginning, solutions presented above seem very promising, in terms of both payload mass and constraint enforcement.

The exploitation of an effective multi-population meta-heuristic optimization code allows a wide exploration of the solution space, and an optimal solution is found in a relatively easy and fast manner. Moreover, handling constraints with the proposed technique does not severely affect the search for the optimal solution, by letting the code explore the entire solution space in order to find a promising feasible region where to move.
Obviously, much more can be done in order to improve the performance of the proposed methodology and increase both the reliability of the solutions and the accuracy of the employed model.
As an example, different steering strategy for both the third stage and the orbital control system can be analyzed, through the use of more realistic guidance laws which resemble the optimal ones. 

The \textit{Vega Light} performance along the obtained ascent trajectory lets one foresee the competitiveness of such a new European light launch vehicle within an increasingly crowded market. The re-use of reliable and proven technologies, the exploitation of already existing facilities and the know-how of the expert engineers involved in \textit{Vega} program, will surely make this new launcher one of the strengths of the upcoming \textit{Vega} family.

\bibliographystyle{AAS_publication}   
\bibliography{references}   

\end{document}